\documentclass[leqno, fleqn]{amsart}%
\usepackage{amsmath}
\usepackage{amsfonts}
\usepackage{amssymb}%

\makeatletter \DeclareMathSymbol{\Gamma}{\mathalpha}{letters}{"00}
\DeclareMathSymbol{\Delta}{\mathalpha}{letters}{"01}
\DeclareMathSymbol{\Theta}{\mathalpha}{letters}{"02}
\DeclareMathSymbol{\Lambda}{\mathalpha}{letters}{"03}
\DeclareMathSymbol{\Xi}{\mathalpha}{letters}{"04}
\DeclareMathSymbol{\Pi}{\mathalpha}{letters}{"05}
\DeclareMathSymbol{\Sigma}{\mathalpha}{letters}{"06}
\DeclareMathSymbol{\Upsilon}{\mathalpha}{letters}{"07}
\DeclareMathSymbol{\Phi}{\mathalpha}{letters}{"08}
\DeclareMathSymbol{\Psi}{\mathalpha}{letters}{"09}
\DeclareMathSymbol{\Omega}{\mathalpha}{letters}{"0A}
\DeclareMathSymbol{\varGamma}{\mathalpha}{operators}{"00}
\DeclareMathSymbol{\varDelta}{\mathalpha}{operators}{"01}
\DeclareMathSymbol{\varTheta}{\mathalpha}{operators}{"02}
\DeclareMathSymbol{\varLambda}{\mathalpha}{operators}{"03}
\DeclareMathSymbol{\varXi}{\mathalpha}{operators}{"04}
\DeclareMathSymbol{\varPi}{\mathalpha}{operators}{"05}
\DeclareMathSymbol{\varSigma}{\mathalpha}{operators}{"06}
\DeclareMathSymbol{\varUpsilon}{\mathalpha}{operators}{"07}
\DeclareMathSymbol{\varPhi}{\mathalpha}{operators}{"08}
\DeclareMathSymbol{\varPsi}{\mathalpha}{operators}{"09}
\DeclareMathSymbol{\varOmega}{\mathalpha}{operators}{"0A}

\newcommand{\allmodesymb}[2]{\relax\ifmmode{\mathchoice
{\mbox{\fontsize{\tf@size}{\tf@size}#1{#2}}}
{\mbox{\fontsize{\tf@size}{\tf@size}#1{#2}}}
{\mbox{\fontsize{\sf@size}{\sf@size}#1{#2}}}
{\mbox{\fontsize{\ssf@size}{\ssf@size}#1{#2}}}} \else \mbox{#1{#2}}\fi}

\makeatother
\makeatletter
\renewcommand*\subjclass[2][2000]{%
  \def\@subjclass{#2}%
  \@ifundefined{subjclassname@#1}{%
    \ClassWarning{\@classname}{Unknown edition (#1) of Mathematics%
      Subject Classification; using '2000'.}%
  }{%
    \@xp\let\@xp\subjclassname\csname subjclassname@#1\endcsname%
  }%
} \makeatother

\theoremstyle{plain}
\newtheorem{theorem}{Theorem}
\newtheorem{proposition}{Proposition}
\newtheorem{lemma}{Lemma}

\newtheorem*{conjecture}{Conjecture}

\theoremstyle{remark}

\newtheorem{remark}{Remark}

\allowdisplaybreaks
\numberwithin{equation}{section}
\begin{document}
\title[K\"{a}hler--Einstein metric]{The K\"{a}hler--Einstein metric for some
Hartogs domains over bounded symmetric domains}
\date{12th January 2005}
\author{An WANG}
\address{A.W.: Dept. of Math., Capital Normal Univ., Beijing 100037, China}
\email{wangancn@sina.com}
\author{Weiping YIN}
\address{W.Y.: Dept. of Math., Capital Normal Univ., Beijing 100037, China}
\email{wyin@mail.cnu.edu.cn}
\author{Liyou ZHANG}
\address{L.Zh.: Dept. of Math., Capital Normal Univ., Beijing 100037, China}
\email{zhangly@mail.cnu.edu.cn}
\author{Guy ROOS}
\address{G.R.: Nevski prospekt 113/4-53, 191024 St Petersburg, Russian Federation}
\email{guy.roos@normalesup.org}
\thanks{Research supported in part by NSF of China (Grant 10471097)}
\thanks{The present work was completed during the stay of G. Roos at Capital Normal
University of Beijing in August-September 2004}

\begin{abstract}
We study the complete K\"{a}hler-Einstein metric of a Hartogs
domain $\widetilde{\Omega}$ built on an irreducible bounded
symmetric domain $\Omega $, using a power $N^{\mu}$ of the generic
norm of $\Omega$. The generating function of the
K\"{a}hler-Einstein metric satisfies a complex Monge-Amp\`{e}re
equation with boundary condition. The domain $\widetilde {\Omega}$
is in general not homogeneous, but it has a subgroup of
automorphisms, the orbits of which are parameterized by
$X\in\lbrack0,1[$. This allows to reduce the Monge-Amp\`{e}re
equation to an ordinary differential equation with limit
condition. This equation can be explicitly solved for a special
value $\mu_{0}$ of $\mu$. We work out the details for the two
exceptional symmetric domains. The special value $\mu_{0}$ seems
also to be significant for the properties of other invariant
metrics like the Bergman metric; a conjecture is stated, which is
proved for the exceptional domains.

\end{abstract}
\maketitle

\section*{Introduction}

Let $D$ be a bounded domain in $\mathbb{C}^{n}$. The \emph{complete
(normalized) K\"{a}hler-Einstein metric} on $D$ is the Hermitian metric $E$%
\[
E_{z}(u,v)=\left.  \partial_{u}\overline{\partial}_{v}g\right\vert _{z},
\]
whose \emph{generating function} $g$ is the unique solution of the complex
Monge-Amp\`{e}re equation with boundary condition%
\begin{align*}
&  \det\left(  \frac{\partial^{2}g}{\partial z^{i}\partial\overline{z}^{j}%
}\right)  =\operatorname{e}^{(n+1)g}\qquad(z\in D),\\
&  g(z) \rightarrow\infty\qquad\qquad(z\rightarrow\partial D)
\end{align*}
(see \cite{ChengYau1980}, \cite{MokYau1983}, \cite{Wu1993}).

Let $\Omega$ be a bounded irreducible symmetric domain in $V\simeq
\mathbb{C}^{d}$; we will always consider such a domain in its \emph{circled}
realization. For a real positive number $\mu$, let $\widetilde{\Omega}$ be the
Hartogs type domain defined by%
\[
\widetilde{\Omega}=\widetilde{\Omega}_{k}(\mu)=\left\{  \left(  z,Z\right)
\in\Omega\times\mathbb{C}^{k}\mid\left\Vert Z\right\Vert ^{2}<N(z,z)^{\mu
}\right\}  ,
\]
where $N(z,z)$ denotes the \emph{generic norm} of $\Omega$ (see
Appendix \ref{GenericNorm}). The Bergman kernel of
$\widetilde{\Omega}$ has been computed in \cite{YinRoos2003} .

\bigskip

In this paper, we study the complete K\"{a}hler-Einstein metric of
$\widetilde{\Omega}$. The domain $\widetilde{\Omega}$ is in
general not homogeneous, but it has a subgroup of automorphisms,
the orbits of which are parameterized by $X\in\lbrack0,1[$. This
allows to reduce the Monge-Amp\`{e}re equation to an ordinary
differential equation, following a method used in
\cite{WangYinZhangZhang2004}  when $\Omega$ belongs to one of the
four series of classical domains. This equation can be explicitly
solved for a special value $\mu_{0}$ of $\mu$. We first study the
case $k=1$, then generalize the results to any integer $k$. Tables
are given in Appendix \ref{Tables}, allowing the reader to apply
the results to each irreducible bounded symmetric domain. In
Section \ref{Exceptional}, we work out some details for the two
exceptional bounded symmetric domains; the construction and the
main properties of these domains are recalled in Appendix
\ref{ExceptionalDomains}. In Section \ref{Conjecture}, starting
from the example of exceptional domains, we state a conjecture
which links the critical exponent $\mu_{0}$ for the
K\"{a}hler-Einstein metric and the properties of the Bergman
kernel of $\widetilde{\Omega}_{k}(\mu)$.

\section{Symmetric domains inflated by discs}

This is the case $k=1$, with
\[
\widetilde{\Omega}=\widetilde{\Omega}_{1}(\mu)=\left\{  \left(  z,w\right)
\in\Omega\times\mathbb{C}\mid\left\vert w\right\vert ^{2}<N(z,z)^{\mu
}\right\}  .
\]

\subsection{Automorphisms}

Let $\Omega$ be a bounded irreducible symmetric domain, $\mu>0$ a real number
and
\[
\widetilde{\Omega}=\widetilde{\Omega}_{1}(\mu)=\left\{  \left(  z,w\right)
\in\Omega\times\mathbb{C}\mid\left\vert w\right\vert ^{2}<N(z,z)^{\mu
}\right\}  .
\]
Let $X$ be the function $X:\widetilde{\Omega}\rightarrow\lbrack0,1[$ defined
by
\begin{equation}
X(z,w)=\frac{\left\vert w\right\vert ^{2}}{N(z,z)^{\mu}}. \label{EQN7}%
\end{equation}
Denote by $\operatorname{Aut}^{\prime}\widetilde{\Omega}$ the subgroup of
automorphisms of $\widetilde{\Omega}$ which leave $X$ invariant.

Let $\Phi\in\operatorname{Aut}\Omega$. Denote by $\operatorname{d}\Phi(z)$ the
differential of $\Phi$ at $z$ and by $J\Phi(z)=\det\operatorname{d}\Phi(z)$
its Jacobian.

\begin{lemma}
\label{L2}Let $\Omega$ be a bounded irreducible symmetric domain in
$V\simeq\mathbb{C}^{d}$, with generic norm $N$ and genus $\gamma$ (see
Appendix \ref{GenericNorm}). Then the function on $\Omega\times\Omega$%
\[
\frac{N(z,z)N(t,t)}{\left\vert N(z,t)\right\vert ^{2}}%
\]
is invariant by each $\Phi\in\operatorname{Aut}\Omega$ acting diagonally on
$\Omega\times\Omega$:%
\begin{equation}
\frac{N(\Phi z,\Phi z)N(\Phi t,\Phi t)}{\left\vert N(\Phi z,\Phi t)\right\vert
^{2}}=\frac{N(z,z)N(t,t)}{\left\vert N(z,t)\right\vert ^{2}}. \label{EQN1}%
\end{equation}

\end{lemma}

\begin{proof}
If $\Phi\in\operatorname{Aut}_{0}\Omega$, the identity component of
$\operatorname{Aut}\Omega$, we deduce from (\ref{K12}):
\[
B(\Phi z,\Phi t)=\operatorname{d}\Phi(z)\circ B(z,t)\circ\operatorname{d}%
\Phi(t)^{\ast}%
\]
(see Appendix \ref{JTS}) and $\det B(z,t)=N(z,t)^{\gamma}$, that
\begin{equation}
N(\Phi z,\Phi t)^{\gamma}=J\Phi(z)N(z,t)^{\gamma}\overline{J\Phi(t)}.
\label{EQN3}%
\end{equation}
If $\Phi\in\operatorname{Aut}\Omega$, it can be written $\Phi=\Phi_{1}\Phi
_{2}$, where $\Phi_{1}\in\operatorname{Aut}_{0}\Omega$ and $\Phi_{2}(0)=0$;
then $\Phi_{2}$ is linear, unitary (with respect to the Bergman metric at $0$)
and leaves $N$ invariant, so that (\ref{EQN3}) holds for $\Phi_{2}$, hence
also for every $\Phi\in\operatorname{Aut}\Omega$. The relation (\ref{EQN1})
follows immediately.
\end{proof}

\begin{proposition}
\label{P1}The group $\operatorname{Aut}^{\prime}\widetilde{\Omega}$ consists
of all $\Psi=\left(  \Psi_{1},\Psi_{2}\right)  $:%
\begin{align}
\Psi_{1}(z,w)  &  =\Phi(z),\label{EQN55}\\
\Psi_{2}(z,w)  &  =w\psi(z), \label{EQN56}%
\end{align}
such that $\Phi\in\operatorname{Aut}\Omega$ and
\begin{equation}
\left\vert \psi(z)\right\vert ^{2}=\left(  \frac{N(\Phi z,\Phi z)}%
{N(z,z)}\right)  ^{\mu}. \label{EQN2}%
\end{equation}
For $\Phi\in\operatorname{Aut}\Omega$, let $z_{0}=\Phi^{-1}(0)$; then the
functions $\psi$ satisfying (\ref{EQN2}) are the functions
\begin{equation}
\psi(z)=\operatorname{e}^{\operatorname{i}\theta}\frac{N\left(  z_{0}%
,z_{0}\right)  ^{\mu/2}}{N\left(  z,z_{0}\right)  ^{\mu}}. \label{EQN4}%
\end{equation}
The orbits of $\operatorname{Aut}^{\prime}\widetilde{\Omega}$ are the level
sets%
\[
\Sigma_{\lambda}=\left\{  X=\lambda\mid\lambda\in\lbrack0,1[\right\}  .
\]

\end{proposition}

\begin{proof}
It is easily checked that the maps $\Psi$ of the form (\ref{EQN55}%
)-(\ref{EQN56}), satisfying (\ref{EQN2}), are automorphisms of $\widetilde
{\Omega}$, leave $X$ invariant and form a subgroup $G$ of automorphisms of
$\widetilde{\Omega}$.

Let $z^{\ast}=\Phi(z)$. Applying (\ref{EQN1}) to $z_{0}=\Phi^{-1}(0)$, we get
\begin{equation}
\frac{N(z^{\ast},z^{\ast})}{N(z,z)}=\frac{N\left(  z_{0},z_{0}\right)
}{\left\vert N(z,z_{0})\right\vert ^{2}}, \label{EQN5}%
\end{equation}
as $N(x,0)=1$ for each $x$. The relation (\ref{EQN3}) also implies
\[
1=J\Phi(z)N\left(  z,z_{0}\right)  ^{\gamma}\overline{J\Phi\left(
z_{0}\right)  };
\]
this means that the holomorphic function $z\mapsto N\left(  z,z_{0}\right)  $
never vanishes on the convex domain $\Omega$, and we can define the
holomorphic function $N\left(  z,z_{0}\right)  ^{\mu}$ for any real $\mu$
assuming it is positive for $z=z_{0}$. By (\ref{EQN5}), the function
\[
\psi_{0}(z)=\frac{N\left(  z_{0},z_{0}\right)  ^{\mu/2}}{N\left(
z,z_{0}\right)  ^{\mu}}%
\]
satisfies (\ref{EQN2}). If $\psi$ is another function satisfying (\ref{EQN2}),
then $\psi/\psi_{0}$ has constant modulus $1$ and $\psi=\operatorname{e}%
^{\operatorname{i}\theta}\psi_{0}$.

Let $\Psi\in\operatorname{Aut}^{\prime}\widetilde{\Omega}$; as $\Psi$
preserves $X$, $\Psi(\Omega\times\{0\})=\Omega\times\{0\}$ and $\Psi
(z,0)=\left(  \Phi(z),0\right)  $ with $\Phi\in\operatorname{Aut}\Omega$.
There exists $\Psi^{1}\in G$ such that $\Psi^{1}(z,0)=\left(  \Phi
(z),0\right)  $ and $\Theta=\Psi\circ\left(  \Psi^{1}\right)  ^{-1}$ is an
element of $\operatorname{Aut}^{\prime}\widetilde{\Omega}$, such that
$\Theta(z,0)=(z,0)$ for all $z\in\Omega$. In particular, $\Theta(0,0)=(0,0)$;
as $\widetilde{\Omega}$ is bounded and circled, it follows from a lemma of
H.~Cartan (see \cite{Cartan1931}, \cite{Loos1977}) that $\Theta$ is linear.
Then $\Theta$ has the form
\begin{align*}
\Theta_{1}(z,w)  &  =z+wu,\\
\Theta_{2}(z,w)  &  =cw.
\end{align*}
The invariance of $X$ under $\Theta$ implies
\[
\left\vert c\right\vert ^{2}=\frac{N(z+wu,z+wu)^{\mu}}{N(z,z)^{\mu}}%
\]
for all $\left(  z,w\right)  \in\widetilde{\Omega}$ and in particular
\[
\left\vert c\right\vert ^{2}=N\left(  wu,wu\right)  ^{\mu}%
\]
for all $w$ such that $\left\vert w\right\vert <1$. This implies $\left\vert
c\right\vert =1$ and $u=0$. Then $\Theta\in G$ and $\Psi=\Theta\circ\Psi^{1}$
belongs also to $G$.

If $\left(  z_{0},w_{0}\right)  $ and $\left(  z_{0}^{\prime},w_{0}^{\prime
}\right)  $ belong to the same level set $\Sigma_{\lambda}$, that is
$X(z_{0},w_{0})=X(z_{0}^{\prime},w_{0}^{\prime})=\lambda$, take $\Phi
\in\operatorname{Aut}\Omega$ such that $\Phi(z_{0})=z_{0}^{\prime}$ and
$\psi:\Omega\rightarrow\mathbb{C}$ such that
\[
\left\vert \psi(z)\right\vert ^{2}=\left(  \frac{N(\Phi z,\Phi z)}%
{N(z,z)}\right)  ^{\mu}.
\]
Then $\Psi\in\operatorname{Aut}^{\prime}\widetilde{\Omega}$ defined by
\begin{align*}
\Psi_{1}(z,w)  &  =\Phi(z),\\
\Psi_{2}(z,w)  &  =w\psi(z)
\end{align*}
maps $\left(  z_{0},w_{0}\right)  $ to $\left(  z_{0}^{\prime},w_{0}\psi
(z_{0})\right)  $; we have $\left\vert w_{0}\psi(z_{0})\right\vert
^{2}=\left\vert w_{0}^{\prime}\right\vert ^{2}$, so it suffices to change
$\Psi_{2}$ in $\alpha\Psi_{2}$ for some $\alpha\in\mathbb{C}$, $\left\vert
\alpha\right\vert =1$, in order to obtain $\Psi(z_{0},w_{0})=\left(
z_{0}^{\prime},w_{0}^{\prime}\right)  $. So the group $\operatorname{Aut}%
^{\prime}\widetilde{\Omega}$ acts transitively on the level sets
$\Sigma_{\lambda}$ of $X$.
\end{proof}

The orbits $\Sigma_{\lambda}$ are real hypersurfaces of $\widetilde{\Omega}$
when $\lambda>0$; the orbit $\Sigma_{0}$ is $\Omega\times\{0\}$. Note also
that $\Psi\in\operatorname{Aut}^{\prime}\widetilde{\Omega}$ extends
continuously to the boundary $\partial\widetilde{\Omega}$ of $\widetilde
{\Omega}$, as $\Phi$ extends continuously to $\partial\Omega$ and
\[
\partial\widetilde{\Omega}=\left(  \partial\Omega\times\{0\}\right)
\cup\left\{  \left(  z,w\right)  \in\Omega\times\mathbb{C}\mid\left\vert
w\right\vert ^{2}=N(z,z)^{\mu}\right\}  .
\]
The part
\[
\partial_{0}\widetilde{\Omega}=\left\{  \left(  z,w\right)  \in\Omega
\times\mathbb{C}\mid\left\vert w\right\vert ^{2}=N(z,z)^{\mu}\right\}
\]
is an orbit of $\operatorname{Aut}^{\prime}\widetilde{\Omega}$ (corresponding
to $X=1$), and $\partial\Omega\times\{0\}$ is a finite union of orbits.

In order to compute the differential of $\psi$ and $\Psi$, we need the
following general result in Jordan triple systems. See Appendix \ref{JTS} for
the notations. In particular, $B$ is the Bergman operator, $N$ the generic
norm, $m_{1}$ the generic trace, $y^{x}$ is the quasi-inverse; if $V$ is
simple (if $\Omega$ is irreducible), the genus is denoted by $\gamma$.

\begin{lemma}
\label{L4}1)~Let $V$ be a Hermitian Jordan triple system. Then
\begin{align}
&  \operatorname{d}_{x}B(x,y) =-D(\operatorname{d}x,y^{x})B(x,y),
\label{EQN8}\\
&  \operatorname{d}_{y}B(x,y) =-B(x,y)D(x^{y},\operatorname{d}y),
\label{EQN8A}\\
&  \frac{\operatorname{d}_{x}\det B(x,y)}{\det B(x,y)} =-\operatorname{tr}%
D(\operatorname{d}x,y^{x}),\label{EQN9}\\
&  \frac{\operatorname{d}_{y}\det B(x,y)}{\det B(x,y)} =-\operatorname{tr}%
D(x^{y},\operatorname{d}y). \label{EQN9A}%
\end{align}

2)~If $V$ is simple,
\begin{align}
&  \frac{\operatorname{d}_{x}N(x,y)}{N(x,y)} =-m_{1}(\operatorname{d}%
x,y^{x}),\label{EQN10}\\
&  \frac{\operatorname{d}_{y}N(x,y)}{N(x,y)} =-m_{1}(x^{y},\operatorname{d}%
y),\label{EQN10A}\\
&  \partial\left(  \frac{\overline{\partial}N(z,z)}{N(z,z)}\right)
=-m_{1}\left(  B(z,z)^{-1}\operatorname{d}z,\operatorname{d}\overline
{z}\right)  . \label{EQN19}%
\end{align}
In particular,
\begin{equation}
\left.  -\partial\left(  \frac{\overline{\partial}N(z,z)}{N(z,z)}\right)
\right\vert _{z=0}=m_{1}\left(  \operatorname{d}z,\operatorname{d}\overline
{z}\right)  . \label{EQN23}%
\end{equation}

\end{lemma}

\begin{proof}
We start from the addition formula for the Bergman operator%
\[
B(x+z,y)=B(z,y^{x})B(x,y)
\]
(see \cite{Roos1999}, p. 469, (J6.4')). Using the definition of $B$, this
identity can be written%
\[
B(x+z,y)=\left(  \operatorname{id}_{V}-D(z,y^{x})+Q(z)Q(y^{x})\right)
B(x,y).
\]
Taking the linear part in $z$ proves (\ref{EQN8}).

In $GL(V)$ we have the well-known relation
\[
\frac{\operatorname{d}(\det A)}{\det A}=\operatorname{tr}\left(
A^{-1}\operatorname{d}A\right)  =\operatorname{tr}\left(  \operatorname{d}%
A.A^{-1}\right)  .
\]
This implies, using (\ref{EQN8}),
\[
\frac{\operatorname{d}_{x}\det B(x,y)}{\det B(x,y)}=\operatorname{tr}\left(
\operatorname{d}_{x}B(x,y).B(x,y)^{-1}\right)  =-\operatorname{tr}%
D(\operatorname{d}x,y^{x}),
\]
which proves (\ref{EQN9}). The \textquotedblleft dual\textquotedblright%
\ formulas (\ref{EQN8A}) and (\ref{EQN9A}) are proved in the same way.

If $V$ is simple, $\det B(x,y)=N(x,y)^{\gamma}$ and $\operatorname{tr}%
D(x,y)=\gamma m_{1}(x,y)$; so (\ref{EQN10})-(\ref{EQN10A}) immediately follow
from (\ref{EQN9})-(\ref{EQN9A}).

As $N(x,y)$ is holomorphic in $x$ and anti-holomorphic in $y$, we have
\[
\frac{\overline{\partial}N(z,z)}{N(z,z)}=-m_{1}(z^{z},\operatorname{d}z)
\]
and
\[
\partial\left(  \frac{\overline{\partial}N(z,z)}{N(z,z)}\right)
=-m_{1}(\partial\left(  z^{z}\right)  ,\operatorname{d}z).
\]
The differential of the quasi-inverse $x^{y}$ with respect to $x$ is (see
\cite{Roos1999}, p. 471, relation (D2))
\[
\operatorname{d}_{x}\left(  x^{y}\right)  =B(x,y)^{-1}\operatorname{d}x;
\]
as $x^{y}$ is holomorphic in $x$ and anti-holomorphic in $y$, this implies%
\[
\partial\left(  z^{z}\right)  =B(z,z)^{-1}\operatorname{d}z
\]
and (\ref{EQN19}).
\end{proof}

\begin{lemma}
\label{L3}1)~Let $\Phi\in\operatorname{Aut}_{0}\Omega$, $\Phi(z_{0})=0$ and
let $\psi:\Omega\rightarrow\mathbb{C}$ defined by (\ref{EQN4}):%
\[
\psi(z)=\operatorname{e}^{\operatorname{i}\theta}\frac{N\left(  z_{0}%
,z_{0}\right)  ^{\mu/2}}{N\left(  z,z_{0}\right)  ^{\mu}}.
\]
Then
\begin{equation}
\operatorname{d}\psi(z)=\mu\psi(z)m_{1}\left(  \operatorname{d}z,z_{0}%
^{\ z}\right)  . \label{EQN11}%
\end{equation}

2)~Let $\Psi=\left(  \Psi_{1},\Psi_{2}\right)  \in\operatorname{Aut}^{\prime
}\widetilde{\Omega}$ be defined by%
\begin{align*}
\Psi_{1}(z,w)  &  =\Phi(z),\\
\Psi_{2}(z,w)  &  =w\psi(z).
\end{align*}
Then the differential of $\Psi$ is given by
\begin{align}
\operatorname{d}\Psi_{1}(z,w)  &  =\operatorname{d}\Phi(z),\label{EQN13}\\
\operatorname{d}\Psi_{2}(z,w)  &  =\mu w\psi(z)m_{1}\left(  \operatorname{d}%
z,z_{0}^{\ z}\right)  +\psi(z)\operatorname{d}w. \label{EQN14}%
\end{align}
The Jacobian of $\Psi$ is
\begin{equation}
J\Psi(z,w)=\psi(z)J\Phi(z) \label{EQN15}%
\end{equation}
and satisfies
\begin{equation}
\left\vert J\Psi(z,w)\right\vert ^{2}=\left(  \frac{N(\Phi z,\Phi z)}%
{N(z,z)}\right)  ^{\gamma+\mu}. \label{EQN16}%
\end{equation}

\end{lemma}

The relation (\ref{EQN11}) follows immediately from (\ref{EQN10}), and implies
(\ref{EQN14}). The triangular form of (\ref{EQN13})-(\ref{EQN14}) yields
(\ref{EQN15}). The relation (\ref{EQN16}) follows then from (\ref{EQN3}) (for
$z=t$) and (\ref{EQN2}).

\begin{remark}
The relation (\ref{EQN19}) expresses the well-known fact that
\[
-\partial\left(  \frac{\overline{\partial}N(z,z)}{N(z,z)}\right)
\]
is equal, up to the factor $\frac{1}{\gamma}$, to the Bergman metric of
$\Omega$ at $z$. This could also be established directly, as the Bergman
kernel of $\Omega$ is
\[
\mathcal{K}(z)=\frac{1}{\operatorname{vol}\Omega}N(z,z)^{-\gamma}.
\]
In this paper, we will only make use of the special case (\ref{EQN23}) at
$z=0$.
\end{remark}

\subsection{Reduction of the Monge-Amp\`{e}re equation}

\subsubsection{}

Let $\Omega$ be a bounded irreducible symmetric domain, $\mu>0$ a real number
and
\[
\widetilde{\Omega}=\widetilde{\Omega}_{1}(\mu)=\left\{  \left(  z,w\right)
\in\Omega\times\mathbb{C}\mid\left\vert w\right\vert ^{2}<N(z,z)^{\mu
}\right\}  .
\]
We denote by $d$ the complex dimension of $\Omega$ and by $n=d+1$ the
dimension of $\widetilde{\Omega}$.

Let $g$ be a $C^{2}$ function in $\widetilde{\Omega}$, which is a solution of
the Monge-Amp\`{e}re equation
\begin{equation}
\det\left(  \frac{\partial^{2}g}{\partial z^{i}\partial\overline{z}^{j}%
}\right)  =\operatorname{e}^{(n+1)g} \label{EQN17}%
\end{equation}
and which generates an invariant form $\partial\overline{\partial}g$. For
$\Psi\in\operatorname{Aut}^{\prime}\widetilde{\Omega}$, the invariance of the
metric and the Monge-Amp\`{e}re equation (\ref{EQN17}) imply
\[
\operatorname{e}^{(n+1)g(z,w)}=\left\vert J\Psi(z,w)\right\vert ^{2}%
\operatorname{e}^{(n+1)g(\Psi(z,w))}%
\]
and, using (\ref{EQN16}),
\[
\operatorname{e}^{(n+1)g(z,w)}N(z,z)^{\gamma+\mu}=\operatorname{e}%
^{(n+1)g(\Psi(z,w))}N(\Phi z,\Phi z)^{\gamma+\mu}.
\]
In other words, \emph{the function}
\begin{equation}
g(z,w)+\frac{\gamma+\mu}{n+1}\log N(z,z) \label{EQN18}%
\end{equation}
\emph{is constant on the orbits of }$\operatorname{Aut}^{\prime}%
\widetilde{\Omega}$. Recall that these orbits are also
parameterized by $X\in\lbrack0,1[$ and \emph{define }$h(X)$ as the
value of the function (\ref{EQN18}) on the orbit $X=\left\vert
w\right\vert ^{2}N(z,z)^{-\mu}$. The
function $g$ can then be written%
\begin{equation}
g(z,w)=-\frac{\gamma+\mu}{n+1}\log N(z,z)+h\left(  \frac{\left\vert
w\right\vert ^{2}}{N(z,z)^{\mu}}\right)  ; \label{EQN20}%
\end{equation}
as $N(0,0)=1$, the function $h$ can be obtained from $g$ by
\[
h\left(  \left\vert w\right\vert ^{2}\right)  =g(0,w).
\]

We will show that the Monge-Amp\`{e}re equation (\ref{EQN17}) is equivalent to
an ordinary differential equation for the function $h$.

\subsubsection{}

Let $\left(  z^{1},\ldots,z^{d},z^{d+1}=w\right)  $ be linear coordinates on
$\Omega\times\mathbb{C}$ and let
\[
\omega=\operatorname{d}z^{1}\wedge\operatorname{d}\overline{z}^{1}\wedge
\cdots\wedge\operatorname{d}z^{d}\wedge\operatorname{d}\overline{z}^{d}%
\wedge\operatorname{d}w\wedge\operatorname{d}\overline{w}.
\]
\emph{We choose for }$\left(  z^{1},\ldots,z^{d}\right)  $\emph{ orthonormal
coordinates w.r. to the Hermitian metric }$m_{1}$\emph{ relative to }$\Omega
$\emph{.} By Lemma \ref{L4}, we have
\[
\left.  -\partial\left(  \frac{\overline{\partial}N(z,z)}{N(z,z)}\right)
\right\vert _{z=0}=m_{1}\left(  \operatorname{d}z,\operatorname{d}\overline
{z}\right)
\]
and%
\begin{equation}
\left.  \frac{\left(  -\partial\overline{\partial}\log N(z,z)\right)  ^{d}%
}{d!}\right\vert _{z=0}=\operatorname{d}z^{1}\wedge\operatorname{d}%
\overline{z}^{1}\wedge\cdots\wedge\operatorname{d}z^{d}\wedge\operatorname{d}%
\overline{z}^{d}. \label{EQN24}%
\end{equation}
For any $C^{2}$ function $f$ on $\widetilde{\Omega}$, we have
\[
\frac{1}{(d+1)!}\left(  \partial\overline{\partial}f\right)  ^{d+1}%
=\det\left(  \frac{\partial^{2}f}{\partial z^{i}\partial\overline{z}^{j}%
}\right)  \omega.
\]

\begin{lemma}
\label{L1}Let $g:\widetilde{\Omega}\rightarrow\mathbb{R}$ be a $C^{2}$
function and let $h:[0,1[\rightarrow\mathbb{R}$ be related to $g$ by
(\ref{EQN20}). Then, for $X=\left\vert w\right\vert ^{2}>0$,
\begin{equation}
\frac{\left(  \partial\overline{\partial}g\right)  ^{d+1}}{(d+1)!}%
(0,w)=\left(  \mu Xh^{\prime}(X)+\frac{\gamma+\mu}{d+2}\right)  ^{d}\left(
Xh^{\prime}(X)\right)  ^{\prime}\omega. \label{EQN22}%
\end{equation}

\end{lemma}

\begin{proof}
Let $X=\left\vert w\right\vert ^{2}N(z,z)^{-\mu}$. Then
\begin{align*}
\frac{\partial X}{X}  &  =\frac{\operatorname{d}w}{w}-\mu\frac{\partial N}%
{N},\\
\frac{\overline{\partial}X}{X}  &  =\frac{\operatorname{d}\overline{w}%
}{\overline{w}}-\mu\frac{\overline{\partial}N}{N}.
\end{align*}
We have
\[
\overline{\partial}\left(  h(X)\right)  =h^{\prime}(X)\overline{\partial
}X=Xh^{\prime}(X)\left(  \frac{\operatorname{d}\overline{w}}{\overline{w}}%
-\mu\frac{\overline{\partial}N}{N}\right)
\]
and
\begin{align*}
\partial\overline{\partial}\left(  h(X)\right)  =  &  \left(  Xh^{\prime
}(X)\right)  ^{\prime}X\,\left(  \frac{\operatorname{d}w}{w}-\mu\frac{\partial
N}{N}\right)  \wedge\left(  \frac{\operatorname{d}\overline{w}}{\overline{w}%
}-\mu\frac{\overline{\partial}N}{N}\right) \\
&  -\mu Xh^{\prime}(X)\partial\left(  \frac{\overline{\partial}N}{N}\right)  .
\end{align*}
We conclude that
\begin{align*}
\partial\overline{\partial}g=  &  -\left(  \mu Xh^{\prime}(X)+\frac{\gamma
+\mu}{d+2}\right)  \partial\overline{\partial}\log N(z,z)\\
&  +\left(  Xh^{\prime}(X)\right)  ^{\prime}X\,\left(  \frac{\operatorname{d}%
w}{w}-\mu\frac{\partial N}{N}\right)  \wedge\left(  \frac{\operatorname{d}%
\overline{w}}{\overline{w}}-\mu\frac{\overline{\partial}N}{N}\right)
\end{align*}
and
\begin{align*}
\left(  \partial\overline{\partial}g\right)  ^{d+1}=(  &  d+1)\left(  \mu
Xh^{\prime}(X)+\frac{\gamma+\mu}{d+2}\right)  ^{d}\left(  Xh^{\prime
}(X)\right)  ^{\prime}\left(  -\partial\overline{\partial}\log N(z,z)\right)
^{d}\\
&  \wedge\frac{X}{\left\vert w\right\vert ^{2}}\operatorname{d}w\wedge
\operatorname{d}\overline{w}.
\end{align*}
At $z=0$, we have $X=\left\vert w\right\vert ^{2}$ and (\ref{EQN24}):
\[
\left.  \frac{\left(  -\partial\overline{\partial}\log N(z,z)\right)  ^{d}%
}{d!}\right\vert _{z=0}=\operatorname{d}z^{1}\wedge\operatorname{d}%
\overline{z}^{1}\wedge\cdots\wedge\operatorname{d}z^{d}\wedge\operatorname{d}%
\overline{z}^{d},
\]
which implies the result (\ref{EQN22}).
\end{proof}

\begin{lemma}
\label{L5}Let $g$ be a $C^{2}$ function on $\widetilde{\Omega}$, which is a
solution of the Monge-Amp\`{e}re equation
\[
\det\left(  \frac{\partial^{2}g}{\partial z^{i}\partial\overline{z}^{j}%
}\right)  =\operatorname{e}^{(n+1)g}%
\]
and which generates an invariant form $\partial\overline{\partial}g$. Let
\begin{align*}
g(z,w)  &  =-\frac{\gamma+\mu}{d+2}\log N(z,z)+h\left(  X\right)  ,\\
X  &  =\frac{\left\vert w\right\vert ^{2}}{N(z,z)^{\mu}}.
\end{align*}
Then $h$ satisfies on $]0,1[$ the differential equation
\begin{align}
&  \left(  \mu Xh^{\prime}(X)+\frac{\gamma+\mu}{d+2}\right)  ^{d}\left(
Xh^{\prime}(X)\right)  ^{\prime}=\operatorname{e}^{(d+2)h},\label{EQN25}\\
&  Xh^{\prime}(X)\rightarrow0\qquad(X\rightarrow0). \label{EQN25A}%
\end{align}
The boundary condition%
\[
g(z)\rightarrow\infty\qquad(z\rightarrow\partial\widetilde{\Omega})
\]
implies
\begin{equation}
h(X)\rightarrow\infty\qquad(X\rightarrow1). \label{EQN26}%
\end{equation}

\end{lemma}

\begin{proof}
The differential equation (\ref{EQN25}) results directly from (\ref{EQN22}).
The limit condition (\ref{EQN26}) results from the boundary condition on $g$,
as
\[
h\left(  \left\vert w\right\vert ^{2}\right)  =g(0,w)
\]
and $(0,1)\in\partial\widetilde{\Omega}$. From
\[
\overline{w}h^{\prime}\left(  \left\vert w\right\vert ^{2}\right)
=\frac{\partial g}{\partial w}(0,w),
\]
we deduce that $X^{1/2}h^{\prime}(X)\rightarrow0$ as $X\rightarrow0.$
\end{proof}

\subsubsection{}

Let
\begin{align}
\beta &  =\frac{\gamma+\mu}{\mu\left(  d+2\right)  },\label{EQN33}\\
Y  &  =Xh^{\prime}(X)+\beta. \label{EQN31}%
\end{align}
It results from the previous lemmas that if $g$ is a solution of the
Monge-Amp\`{e}re equation%
\[
\det\left(  \frac{\partial^{2}g}{\partial z^{i}\partial\overline{z}^{j}%
}\right)  =\operatorname{e}^{(d+2)g}%
\]
with the boundary condition
\[
g(z)\rightarrow\infty\qquad(z\rightarrow\partial\widetilde{\Omega}),
\]
the function $Y$ defined by (\ref{EQN31}) satisfies the differential equation
\begin{equation}
\left(  \mu Y\right)  ^{d}Y^{\prime}=\operatorname{e}^{(d+2)h} \label{EQN32}%
\end{equation}
with the initial condition $Y(0)=\beta$.

This shows that the derivative of $Y^{d+1}$ is positive and tends to $\infty$
when $X\rightarrow1$, as $h(X)\rightarrow\infty$ when $X\rightarrow1$. So the
function $Y$ is strictly increasing and maps $[0,1[$ onto $[\beta
,+\infty\lbrack$.

Taking logarithmic derivatives of both sides of (\ref{EQN32}), we get
\begin{equation}
\frac{\left(  Y^{d}Y^{\prime}\right)  ^{\prime}}{Y^{d}Y^{\prime}%
}=(d+2)h^{\prime} \label{EQN32B}%
\end{equation}
and, using the definition of $Y$,
\begin{equation}
\frac{\left(  Y^{d}Y^{\prime}\right)  ^{\prime}}{Y^{d}Y^{\prime}}%
=(d+2)\frac{Y-\beta}{X}. \label{EQN32A}%
\end{equation}
This can be written
\[
X\left(  Y^{d}Y^{\prime}\right)  ^{\prime}=(d+2)Y^{d}\left(  Y-\beta\right)
Y^{\prime}%
\]
or
\begin{align*}
\left(  XY^{d}Y^{\prime}\right)  ^{\prime}  &  =Y^{d}Y^{\prime}+(d+2)Y^{d}%
\left(  Y-\beta\right)  Y^{\prime}\\
&  =(d+2)Y^{d+1}Y^{\prime}-\frac{\gamma}{\mu}Y^{d}Y^{\prime},
\end{align*}
using the definition (\ref{EQN33}) of $\beta$. From (\ref{EQN25A}), we deduce
that $Y(0)=\beta$; integrating with this initial condition yields%
\begin{equation}
XY^{d}Y^{\prime}=Y^{d+2}-\beta^{d+2}-\frac{\gamma}{\mu\left(  d+1\right)
}\left(  Y^{d+1}-\beta^{d+1}\right)  . \label{EQN34}%
\end{equation}

Let us denote by $P$ the polynomial
\begin{equation}
P(Y)=Y^{d+2}-\beta^{d+2}-\frac{\gamma}{\mu\left(  d+1\right)  }\left(
Y^{d+1}-\beta^{d+1}\right)  . \label{EQN34A}%
\end{equation}
By construction,
\[
P^{\prime}(Y)=Y^{d}+(d+2)Y^{d}\left(  Y-\beta\right)  ,
\]
which shows that $P^{\prime}(t)>0$ for $t>\beta$ and consequently $P(t)>0$ for
$t>\beta$. Let $R$ be the polynomial defined by
\[
P(Y)=\left(  Y-\beta\right)  R(Y);
\]
then $R(\beta)=\beta^{d}$ and $R$ is strictly positive on $[\beta
,\infty\lbrack$.

\subsubsection{}

Now we prove that the resolution of the ordinary differential equation
(\ref{EQN34}) allows to construct the generating function $g$ for the
K\"{a}hler-Einstein metric of $\widetilde{\Omega}=\widetilde{\Omega}_{1}(\mu)$.

\begin{lemma}
\label{L6}The differential equation
\begin{align}
&  XY^{d}Y^{\prime} =P(Y),\label{EQN36}\\
&  Y \rightarrow\infty\qquad\qquad(X\rightarrow1), \label{EQN37}%
\end{align}
where the polynomial $P$ is defined by (\ref{EQN34A}), has a unique solution
\[
Y:[0,1[\rightarrow\lbrack\beta,+\infty\lbrack.
\]
This solution is $C^{\infty}$ at $0$.
\end{lemma}

\begin{proof}
Let $Y:]c,1[\rightarrow\lbrack\mathbb{\beta},\infty\lbrack$ satisfy
(\ref{EQN36})-(\ref{EQN37}). As $P$ is positive on $]\beta,\infty\lbrack$, the
function $Y$ is monotone and its inverse function satisfies the differential
equation%
\begin{align*}
&  \frac{1}{X}\frac{\operatorname{d}X}{\operatorname{d}Y} =\frac{Y^{d}}%
{P(Y)},\\
&  X \rightarrow1\qquad\qquad(Y\rightarrow\infty).
\end{align*}
The solution of this equation is given by%
\begin{equation}
-\log X={\displaystyle\int\limits_{Y}^{\infty}} \frac{y^{d}\operatorname{d}%
y}{P(y)}. \label{EQN40}%
\end{equation}
This gives $X$ as a function of $Y$, and $Y$ as an implicit function of $X$.
It is defined on $]\mathbb{\beta},\infty\lbrack$ and maps $]\beta
,\infty\lbrack$ on $[0,1[$, as ${\displaystyle\int\limits_{\beta}^{\infty}}
\frac{y^{d}\operatorname{d}y}{P(y)}=+\infty$. So the maximal solution of
(\ref{EQN36})-(\ref{EQN37}) is defined on $]0,1[$; it is $C^{\infty}$ on
$]0,1[$ and extends continuously to $[0,1[$, with $Y(0)=\beta$.

The relation (\ref{EQN40}) can be written%
\begin{equation}
-\log X=-\log(Y-\beta)+\log\beta+{\displaystyle\int\limits_{Y}^{2\beta}}
\frac{\left(  y^{d}-R(y)\right)  \operatorname{d}y}{(y-\beta)R(y)}%
+{\displaystyle\int\limits_{2\beta}^{\infty}} \frac{y^{d}\operatorname{d}%
y}{P(y)}. \label{EQN42A}%
\end{equation}
The polynomial $R$ is positive on $[\beta,\infty\lbrack$ and $R(\beta
)=P^{\prime}(\beta)=\beta^{d}$. Let $S$ be defined by
\[
y^{d}-R(y)=(y-\beta)S(y).
\]
Then (\ref{EQN42A}) may be written, for $Y>\beta$,
\[
-\log X=-\log(Y-\beta)+C_{0}-{\displaystyle\int\limits_{\beta}^{Y}}
\frac{S(y)\operatorname{d}y}{R(y)},
\]
with%
\[
C_{0}=\log\beta+{\displaystyle\int\limits_{\beta}^{2\beta}} \frac{\left(
y^{d}-R(y)\right)  \operatorname{d}y}{(y-\beta)R(y)}+{\displaystyle\int
\limits_{2\beta}^{\infty}} \frac{y^{d}\operatorname{d}y}{P(y)}.
\]
For $Y>\beta$, we have
\[
X=\operatorname{e}^{-C_{0}}(Y-\beta)\exp{\displaystyle\int\limits_{\beta}^{Y}}
\frac{S(y)\operatorname{d}y}{R(y)},
\]
which shows that $X$ is a $C^{\infty}$ invertible function of $Y\in
\lbrack\beta,\infty\lbrack$, and $Y$ a $C^{\infty}$ function of $X\in
\lbrack0,1[$.
\end{proof}

\begin{theorem}
\label{T1}The generating function $g$ for the K\"{a}hler-Einstein metric of
\[
\widetilde{\Omega}=\widetilde{\Omega}(\mu)=\left\{  \left(  z,w\right)
\in\Omega\times\mathbb{C}\mid\left\vert w\right\vert ^{2}<N(z,z)^{\mu
}\right\}
\]
is given by
\begin{equation}
g(z,w)=-\frac{\gamma+\mu}{d+2}\log N(z,z)+h\left(  \frac{\left\vert
w\right\vert ^{2}}{N(z,z)^{\mu}}\right)  , \label{EQN42}%
\end{equation}
where
\begin{equation}
\operatorname{e}^{(d+2)h}=\left(  \mu Y\right)  ^{d}Y^{\prime} \label{EQN36A}%
\end{equation}
and the function $Y:[0,1[\rightarrow\lbrack\beta,+\infty\lbrack$ is the
solution of (\ref{EQN36})-(\ref{EQN37}).
\end{theorem}

\begin{proof}
We have already seen that, if $g$ is the generating function for the
K\"{a}hler-Einstein metric of $\widetilde{\Omega}(\mu)$, the functions $h$ and
$Y$ satisfy the conditions of the theorem.

Let now $Y:[0,1[\rightarrow\lbrack\beta,\infty\lbrack$ be the solution of
(\ref{EQN36})-(\ref{EQN37}).and let $h$ and $g$ be defined from $Y$ by
(\ref{EQN36A}) and (\ref{EQN42}). As $Y$ verifies (\ref{EQN36}), which is
equivalent to (\ref{EQN34}), this implies (\ref{EQN32A}). Comparing with the
logarithmic derivative of (\ref{EQN36A}), we get for $X>0$%
\[
\frac{Y-\beta}{X}=h^{\prime}(X).
\]
This can be written
\[
h^{\prime}(X)=\int_{0}^{1}Y^{\prime}(tX)\operatorname{d}t\qquad(X>0),
\]
which implies
\[
h^{\prime\prime}(X)=\int_{0}^{1}XY^{\prime\prime}(tX)\operatorname{d}%
t\qquad(X>0)
\]
and the existence of $h^{\prime}(0)=Y^{\prime}(0)$ and $h^{\prime\prime
}(0)=\frac{1}{2}Y^{\prime\prime}(0)$ follows. So $h$ is $C^{2}$ on $[0,1[$ and
$g$ is $C^{2}$ on $\widetilde{\Omega}.$

For $X\in\lbrack0,1[$, we have then $Y=Xh^{\prime}(X)+\beta$ and $h$ satisfies
the differential equation (\ref{EQN25})
\[
\left(  \mu Xh^{\prime}(X)+\frac{\gamma+\mu}{d+2}\right)  ^{d}\left(
Xh^{\prime}(X)\right)  ^{\prime}=\operatorname{e}^{(d+2)h}.
\]
Using Lemma \ref{L1}, we obtain
\begin{equation}
\det\left(  \frac{\partial^{2}g}{\partial z^{i}\partial\overline{z}^{j}%
}\right)  (0,w)=\operatorname{e}^{(d+2)h\left(  \left\vert w\right\vert
^{2}\right)  }=\operatorname{e}^{(d+2)g(0,w)}, \label{EQN39A}%
\end{equation}
which means that $g$ satisfies the Monge-Amp\`{e}re equation at the points
$\left(  0,w\right)  $. Let $\left(  z,w\right)  \in\widetilde{\Omega}$ and
$\Psi\in\operatorname{Aut}^{\prime}\widetilde{\Omega}$ such that $\Psi\left(
z,w\right)  =\left(  0,w^{\prime}\right)  $. Then%
\begin{align*}
\det\left(  \frac{\partial^{2}g}{\partial z^{i}\partial\overline{z}^{j}%
}\right)  (z,w)  &  =\det\left(  \frac{\partial^{2}g}{\partial z^{i}%
\partial\overline{z}^{j}}\right)  (0,w^{\prime})\left\vert J\Psi
(z,w)\right\vert ^{2}\\
&  =\det\left(  \frac{\partial^{2}g}{\partial z^{i}\partial\overline{z}^{j}%
}\right)  (0,w^{\prime})\frac{1}{N(z,z)^{\gamma+\mu}},\\
\left\vert w^{\prime}\right\vert ^{2}  &  =\frac{\left\vert w\right\vert ^{2}%
}{N(z,z)^{\mu}},\\
g(z,w)  &  =-\frac{\gamma+\mu}{d+2}\log N(z,z)+g(0,w^{\prime}),
\end{align*}
so that (\ref{EQN39A}) implies that $g$ satisfies the Monge-Amp\`{e}re
equation on $\widetilde{\Omega}.$

It remains to prove that $g(z,w)\rightarrow\infty$ as $\left(  z,w\right)
\rightarrow\left(  z_{0},w_{0}\right)  \in\partial\widetilde{\Omega}$. From
$Y=Xh^{\prime}(X)+\beta$ and $Y\rightarrow\infty$ as $X\rightarrow1$, we see
that $h^{\prime}(X)>0$, $h^{\prime}(X)\rightarrow\infty$ and $h(X)\rightarrow
\infty$ as $X\rightarrow1$. The boundary points $\left(  z_{0},w_{0}\right)  $
of $\widetilde{\Omega}$ are of two different types:

\begin{itemize}
\item $z_{0}\in\Omega$, $\left\vert w_{0}\right\vert ^{2}=N\left(  z_{0}%
,z_{0}\right)  ^{\mu}$. Then $X\left(  z_{0},w_{0}\right)  =1$ and%
\[
g(z,w)=-\frac{\gamma+\mu}{d+2}\log N(z,z)+h\left(  X(z,w)\right)
\rightarrow\infty
\]
as $\left(  z,w\right)  \rightarrow\left(  z_{0},w_{0}\right)  $.

\item $z_{0}\in\partial\Omega$, $w_{0}=0$. In this case, $N(z,z)\rightarrow0$
and $h(X)\geq h(0)$, which shows again that $g(z,w)\rightarrow\infty$ as
$\left(  z,w\right)  \rightarrow\left(  z_{0},w_{0}\right)  $.
\end{itemize}
\end{proof}

\begin{remark}
It is easy to check directly that $\partial\overline{\partial}g$ defines a
K\"{a}hler metric. We have
\begin{align*}
\partial\overline{\partial}g=  &  -\left(  \mu Xh^{\prime}(X)+\frac{\gamma
+\mu}{d+2}\right)  \partial\overline{\partial}\log N(z,z)\\
&  +\left(  Xh^{\prime}(X)\right)  ^{\prime}X\,\left(  \frac{\operatorname{d}%
w}{w}-\mu\frac{\partial N}{N}\right)  \wedge\left(  \frac{\operatorname{d}%
\overline{w}}{\overline{w}}-\mu\frac{\overline{\partial}N}{N}\right) \\
=  &  -\mu Y\partial\overline{\partial}\log N(z,z)+XY^{\prime}\left(
\frac{\operatorname{d}w}{w}-\mu\frac{\partial N}{N}\right)  \wedge\left(
\frac{\operatorname{d}\overline{w}}{\overline{w}}-\mu\frac{\overline{\partial
}N}{N}\right)  .
\end{align*}
The associated Hermitian form is
\begin{align*}
H(\zeta,\omega) =  &  -\mu Y\sum\frac{\partial^{2}}{\partial z^{j}%
\partial\overline{z}^{k}}\log N(z,z)\zeta^{j}\overline{\zeta}^{k}\\
&  +XY^{\prime}\left\vert \frac{\omega}{w}-\frac{\mu}{N(z,z)}\sum
\frac{\partial}{\partial z^{j}}N(z,z)\zeta^{j}\right\vert ^{2}.
\end{align*}
The term
\[
B(\zeta)=-\sum\frac{\partial^{2}}{\partial z^{j}\partial\overline{z}^{k}}\log
N(z,z)\zeta^{j}\overline{\zeta}^{k}%
\]
is, up to a constant factor, the Bergman metric of $\Omega$ at $z$. Hence
$H(\zeta,\omega)\geq0$, as $Y>0$ and $Y^{\prime}>0$. If $H(\zeta,\omega)=0$,
then $B(\zeta)=0$, which implies $\zeta=0$ and then $\omega=0$.
\end{remark}

\subsection{The critical exponent}

If $\mu=\mu_{0}$, we have $C=0$ and (\ref{EQN40}) has a very simple form. We
call%
\[
\mu_{0}=\frac{\gamma}{d+1}%
\]
the \emph{critical exponent }for the bounded symmetric domain $\Omega$.

If $\mu=\mu_{0}$, we have $C=0$ and (\ref{EQN40}) is%
\[
-\log X={\displaystyle\int\limits_{Y}^{\infty}}\frac{\operatorname{d}y}%
{y^{2}-y}=\left.  \log\frac{y-1}{y}\right\vert _{Y}^{\infty}=-\log\frac
{Y-1}{Y},
\]
which gives $X=\frac{Y-1}{Y}$ or
\[
Y=\frac{1}{1-X}.
\]
From $\operatorname{e}^{(d+2)h}=\left(  \mu Y\right)  ^{d}Y^{\prime}$, we
obtain
\[
\operatorname{e}^{(d+2)h}=\left(  \mu_{0}\right)  ^{d}\frac{1}{\left(
1-X\right)  ^{d+2}},
\]
that is
\begin{align*}
h  &  =\frac{d}{d+2}\log\mu_{0}+\log\left(  \frac{1}{1-X}\right)  ,\\
&  =\frac{d}{d+2}\log\mu_{0}+\log\left(  \frac{N(z,z)^{\mu_{0}}}%
{N(z,z)^{\mu_{0}}-\left\vert w\right\vert ^{2}}\right)  .
\end{align*}
Here
\[
\frac{\gamma+\mu_{0}}{d+2}=\mu_{0};
\]
applying (\ref{EQN42}), we have%
\begin{align*}
g(z,w)  &  =-\mu_{0}\log N(z,z)+\frac{d}{d+2}\log\mu_{0}+\log\left(
\frac{N(z,z)^{\mu_{0}}}{N(z,z)^{\mu_{0}}-\left\vert w\right\vert ^{2}}\right)
\\
&  =.\frac{d}{d+2}\log\mu_{0}+\log\left(  \frac{1}{N(z,z)^{\mu_{0}}-\left\vert
w\right\vert ^{2}}\right)  .
\end{align*}
The K\"{a}hler-Einstein metric of $\widetilde{\Omega}_{1}\left(  \mu
_{0}\right)  $ is associated to the K\"{a}hler form%
\[
\partial\overline{\partial}g=-\partial\overline{\partial}\log\left(
N(z,z)^{\mu_{0}}-\left\vert w\right\vert ^{2}\right)  .
\]

\section{Inflation by Hermitian balls}

The above results can be extended to the domain
\[
\widetilde{\Omega}=\widetilde{\Omega}_{k}(\mu)=\left\{  \left(  z,Z\right)
\in\Omega\times\mathbb{C}^{k}\mid\left\Vert Z\right\Vert ^{2}<N(z,z)^{\mu
}\right\}  ,
\]
where $\mathbb{C}^{k}$ is endowed with the standard Hermitian norm%
\[
\left\Vert Z\right\Vert ^{2}=\sum_{j=1}^{k}\left\vert Z^{j}\right\vert ^{2}.
\]
We outline the results, omitting the proofs when they are entirely analogous
to the case $k=1$.

\subsection{Automorphisms}

Let $X$ be the function $X:\widetilde{\Omega}\rightarrow\lbrack0,1[$ defined
by
\[
X(z,Z)=\frac{\left\Vert Z\right\Vert ^{2}}{N(z,z)^{\mu}}.
\]
Denote by $\operatorname{Aut}^{\prime}\widetilde{\Omega}$ the subgroup of
automorphisms of $\widetilde{\Omega}$ which leave $X$ invariant.

\begin{proposition}
\label{P2}The group $\operatorname{Aut}^{\prime}\widetilde{\Omega}$ consists
of all $\Psi=\left(  \Psi_{1},\Psi_{2}\right)  $:%
\begin{align*}
\Psi_{1}(z,Z)  &  =\Phi(z),\\
\Psi_{2}(z,Z)  &  =\psi(z)U(Z),
\end{align*}
where $\Phi\in\operatorname{Aut}\Omega$, $U:\mathbb{C}^{k}\rightarrow
\mathbb{C}^{k}$ is special unitary and $\psi$ satisfies
\[
\left\vert \psi(z)\right\vert ^{2}=\left(  \frac{N(\Phi z,\Phi z)}%
{N(z,z)}\right)  ^{\mu}.
\]
The orbits of $\operatorname{Aut}^{\prime}\widetilde{\Omega}$ are the level
sets%
\[
\Sigma_{\lambda}=\left\{  X=\lambda\mid\lambda\in\lbrack0,1[\right\}  .
\]

\end{proposition}

The construction of the functions $\psi$ is given in Proposition \ref{P1}.

\begin{lemma}
\label{L10}Let $\Psi=\left(  \Psi_{1},\Psi_{2}\right)  \in\operatorname{Aut}%
^{\prime}\widetilde{\Omega}$ be defined as above by%
\begin{align*}
\Psi_{1}(z,Z)  &  =\Phi(z),\\
\Psi_{2}(z,Z)  &  =\psi(z)U(Z),
\end{align*}
where $U:\mathbb{C}^{k}\rightarrow\mathbb{C}^{k}$ is special unitary. The
Jacobian of $\Psi$ is
\[
J\Psi(z,Z)=\psi^{k}(z)J\Phi(z)
\]
and satisfies
\begin{equation}
\left\vert J\Psi(z,Z)\right\vert ^{2}=\left(  \frac{N(\Phi z,\Phi z)}%
{N(z,z)}\right)  ^{\gamma+k\mu}. \label{EQN41}%
\end{equation}

\end{lemma}

\subsection{Reduction of the Monge-Amp\`{e}re equation}

Let $\Omega$ be a bounded irreducible symmetric domain, $\mu>0$ a real number
and
\[
\widetilde{\Omega}=\widetilde{\Omega}_{k}(\mu)=\left\{  \left(  z,Z\right)
\in\Omega\times\mathbb{C}^{k}\mid\left\Vert Z\right\Vert ^{2}<N(z,z)^{\mu
}\right\}  .
\]
We denote by $d$ the complex dimension of $\Omega$ and by $n=d+k$ the
dimension of $\widetilde{\Omega}$.

Let $g$ be a $C^{2}$ function in $\widetilde{\Omega}$, which is a solution of
the Monge-Amp\`{e}re equation
\[
\det\left(  \frac{\partial^{2}g}{\partial z^{i}\partial\overline{z}^{j}%
}\right)  _{1\leq i,j\leq d+k}=\operatorname{e}^{(n+1)g}%
\]
and which generates an invariant form $\partial\overline{\partial}g$. For
$\Psi\in\operatorname{Aut}^{\prime}\widetilde{\Omega}$, the invariance of the
metric and the Monge-Amp\`{e}re equation imply
\[
\operatorname{e}^{(n+1)g(z,Z)}=\left\vert J\Psi(z,Z)\right\vert ^{2}%
\operatorname{e}^{(n+1)g(\Psi(z,Z))}%
\]
and, using (\ref{EQN41}),
\[
\operatorname{e}^{(n+1)g(z,Z)}N(z,z)^{\gamma+k\mu}=\operatorname{e}%
^{(n+1)g(\Psi(z,Z))}N(\Phi z,\Phi z)^{\gamma+k\mu}.
\]
The function
\[
g(z,Z)+\frac{\gamma+k\mu}{d+k+1}\log N(z,z)
\]
is then constant on the orbits of\emph{ }$\operatorname{Aut}^{\prime
}\widetilde{\Omega}$. For $X\in\lbrack0,1[$, we define $h(X)$ as the value of
this function on the orbit $X=\left\Vert Z\right\Vert ^{2}N(z,z)^{-\mu}$. The
function $g$ can then be written%
\begin{equation}
g(z,Z)=-\frac{\gamma+k\mu}{d+k+1}\log N(z,z)+h\left(  \frac{\left\Vert
Z\right\Vert ^{2}}{N(z,z)^{\mu}}\right)  ; \label{EQN44}%
\end{equation}
the function $h$ can be obtained from $g$ by
\[
h\left(  \left\Vert Z\right\Vert ^{2}\right)  =g(0,Z)
\]
or
\[
h(X)=g\left(  0,\left(  X^{1/2},0,\ldots,0\right)  \right)  .
\]

Let $(z^{1},\ldots,z^{d})$\emph{ }be coordinates on $V\supset\Omega$, which
are orthonormal w.r. to the Hermitian metric $m_{1}$ relative to $\Omega$ and
let $\left(  Z^{1},\ldots,Z^{k}\right)  $ be orthonormal coordinates for the
Hermitian space $\mathbb{C}^{k}$. Let
\[
\left(  z^{1},\ldots,z^{d+k}\right)  =\left(  z^{1},\ldots,z^{d},Z^{1}%
,\ldots,Z^{k}\right)
\]
and%
\begin{align*}
\omega(z,Z)  &  =\omega_{d}(z)\wedge\omega_{k}(Z)\\
&  =\operatorname{d}z^{1}\wedge\operatorname{d}\overline{z}^{1}\wedge
\cdots\wedge\operatorname{d}z^{d}\wedge\operatorname{d}\overline{z}^{d}%
\wedge\operatorname{d}Z^{1}\wedge\operatorname{d}\overline{Z}^{1}\wedge
\cdots\wedge\operatorname{d}Z^{k}\wedge\operatorname{d}\overline{Z}^{k}.
\end{align*}
For any $C^{2}$ function $f$ on $\widetilde{\Omega}$, we have
\[
\frac{1}{(d+k)!}\left(  \partial\overline{\partial}f\right)  ^{d+k}%
=\det\left(  \frac{\partial^{2}f}{\partial z^{i}\partial\overline{z}^{j}%
}\right)  \omega.
\]

\begin{lemma}
\label{L8}Let $g:\widetilde{\Omega}\rightarrow\mathbb{R}$ and
$h:[0,1[\rightarrow\mathbb{R}$ be $C^{2}$ functions related by (\ref{EQN44}).
Then, for $\left\Vert Z\right\Vert ^{2}=X>0$,%
\begin{equation}
\frac{\left(  \partial\overline{\partial}g\right)  ^{d+k}}{(d+k)!}%
(0,Z)=\left(  h^{\prime}(X)\right)  ^{k-1}\left(  \mu Xh^{\prime}%
(X)+\frac{\gamma+k\mu}{d+k+1}\right)  ^{d}\left(  Xh^{\prime}(X)\right)
^{\prime}\omega. \label{EQN45}%
\end{equation}

\end{lemma}

\begin{proof}
Let $X=\left\Vert Z\right\Vert ^{2}N(z,z)^{-\mu}$. Then
\begin{align*}
\frac{\partial X}{X}  &  =\frac{\partial\left\Vert Z\right\Vert ^{2}%
}{\left\Vert Z\right\Vert ^{2}}-\mu\frac{\partial N}{N},\\
\frac{\overline{\partial}X}{X}  &  =\frac{\overline{\partial}\left\Vert
Z\right\Vert ^{2}}{\left\Vert Z\right\Vert ^{2}}-\mu\frac{\overline{\partial
}N}{N}.
\end{align*}
We have
\[
\overline{\partial}\left(  h(X)\right)  =h^{\prime}(X)\overline{\partial
}X=Xh^{\prime}(X)\left(  \frac{\overline{\partial}\left\Vert Z\right\Vert
^{2}}{\left\Vert Z\right\Vert ^{2}}-\mu\frac{\overline{\partial}N}{N}\right)
\]
and
\begin{align}
\partial\overline{\partial}\left(  h(X)\right)  =  &  \left(  Xh^{\prime
}(X)\right)  ^{\prime}X\,\left(  \frac{\partial\left\Vert Z\right\Vert ^{2}%
}{\left\Vert Z\right\Vert ^{2}}-\mu\frac{\partial N}{N}\right)  \wedge\left(
\frac{\overline{\partial}\left\Vert Z\right\Vert ^{2}}{\left\Vert Z\right\Vert
^{2}}-\mu\frac{\overline{\partial}N}{N}\right) \label{EQN46}\\
&  +Xh^{\prime}(X)\partial\left(  \frac{\overline{\partial}\left\Vert
Z\right\Vert ^{2}}{\left\Vert Z\right\Vert ^{2}}\right)  -\mu Xh^{\prime
}(X)\partial\left(  \frac{\overline{\partial}N}{N}\right)  .\nonumber
\end{align}
As
\[
g(z,Z)=-\frac{\gamma+k\mu}{d+k+1}\log N(z,z)+h\left(  X\right)  ;
\]
we have%
\begin{align*}
\partial\overline{\partial}g=  &  Xh^{\prime}(X)\partial\left(  \frac
{\overline{\partial}\left\Vert Z\right\Vert ^{2}}{\left\Vert Z\right\Vert
^{2}}\right)  -\left(  \mu Xh^{\prime}(X)+\frac{\gamma+k\mu}{d+k+1}\right)
\partial\overline{\partial}\log N(z,z)\\
&  +XY^{\prime}\,\left(  \frac{\partial\left\Vert Z\right\Vert ^{2}%
}{\left\Vert Z\right\Vert ^{2}}-\mu\frac{\partial N}{N}\right)  \wedge\left(
\frac{\overline{\partial}\left\Vert Z\right\Vert ^{2}}{\left\Vert Z\right\Vert
^{2}}-\mu\frac{\overline{\partial}N}{N}\right)  .
\end{align*}
Let
\begin{align*}
\beta &  =\frac{\gamma+k\mu}{\mu\left(  d+k+1\right)  },\\
Y_{0}  &  =Xh^{\prime}(X),\\
Y  &  =Xh^{\prime}(X)+\frac{\gamma+k\mu}{\mu\left(  d+k+1\right)  }%
=Y_{0}+\beta.
\end{align*}
Then
\begin{align*}
\partial\overline{\partial}g =\  &  Y_{0} \partial\left(  \frac{\overline
{\partial}\left\Vert Z\right\Vert ^{2}}{\left\Vert Z\right\Vert ^{2}}\right)
-\mu Y\partial\overline{\partial}\log N(z,z)\\
&  +X\,Y^{\prime}\left(  \frac{\partial\left\Vert Z\right\Vert ^{2}%
}{\left\Vert Z\right\Vert ^{2}}-\mu\frac{\partial N}{N}\right)  \wedge\left(
\frac{\overline{\partial}\left\Vert Z\right\Vert ^{2}}{\left\Vert Z\right\Vert
^{2}}-\mu\frac{\overline{\partial}N}{N}\right)
\end{align*}
and
\begin{align*}
\left(  \partial\overline{\partial}g\right)  ^{d+k}=\  &  \left(
Y_{0}\partial\left(  \frac{\overline{\partial}\left\Vert Z\right\Vert ^{2}%
}{\left\Vert Z\right\Vert ^{2}}\right)  -\mu Y\partial\overline{\partial}\log
N(z,z)\right)  ^{d+k}\\
&  +(d+k)\left(  Y_{0}\partial\left(  \frac{\overline{\partial}\left\Vert
Z\right\Vert ^{2}}{\left\Vert Z\right\Vert ^{2}}\right)  -\mu Y\partial
\overline{\partial}\log N(z,z)\right)  ^{d+k-1}\\
&  \wedge XY^{\prime}\,\left(  \frac{\partial\left\Vert Z\right\Vert ^{2}%
}{\left\Vert Z\right\Vert ^{2}}-\mu\frac{\partial N}{N}\right)  \wedge\left(
\frac{\overline{\partial}\left\Vert Z\right\Vert ^{2}}{\left\Vert Z\right\Vert
^{2}}-\mu\frac{\overline{\partial}N}{N}\right) \\
=\  &  \binom{d+k}{k}\mu^{d}Y_{0}^{k}Y^{d}\left(  -\partial\overline{\partial
}\log N(z,z)\right)  ^{d}\wedge\left(  \partial\left(  \frac{\overline
{\partial}\left\Vert Z\right\Vert ^{2}}{\left\Vert Z\right\Vert ^{2}}\right)
\right)  ^{k}\\
&  +(d+k)\binom{d+k-1}{k-1}\mu^{d}XY^{\prime}Y_{0}^{k-1}Y^{d}\\
&  \left(  -\partial\overline{\partial}\log N(z,z)\right)  ^{d}\wedge\left(
\partial\left(  \frac{\overline{\partial}\left\Vert Z\right\Vert ^{2}%
}{\left\Vert Z\right\Vert ^{2}}\right)  \right)  ^{k-1}\wedge\frac
{\partial\left\Vert Z\right\Vert ^{2}}{\left\Vert Z\right\Vert ^{2}}%
\wedge\frac{\overline{\partial}\left\Vert Z\right\Vert ^{2}}{\left\Vert
Z\right\Vert ^{2}}\\
&  +(d+k)\binom{d+k-1}{k}\mu^{d+2}XY^{\prime}Y_{0}^{k}Y^{d-1}\\
&  \left(  -\partial\overline{\partial}\log N(z,z)\right)  ^{d-1}\wedge
\,\frac{\partial N}{N}\wedge\frac{\overline{\partial}N}{N}\wedge\left(
\partial\left(  \frac{\overline{\partial}\left\Vert Z\right\Vert ^{2}%
}{\left\Vert Z\right\Vert ^{2}}\right)  \right)  ^{k}.
\end{align*}
One checks easily that
\[
\left(  \partial\left(  \frac{\overline{\partial}\left\Vert Z\right\Vert ^{2}%
}{\left\Vert Z\right\Vert ^{2}}\right)  \right)  ^{k}=0
\]
and
\[
\left(  \partial\left(  \frac{\overline{\partial}\left\Vert Z\right\Vert ^{2}%
}{\left\Vert Z\right\Vert ^{2}}\right)  \right)  ^{k-1}\wedge\frac
{\partial\left\Vert Z\right\Vert ^{2}}{\left\Vert Z\right\Vert ^{2}}%
\wedge\frac{\overline{\partial}\left\Vert Z\right\Vert ^{2}}{\left\Vert
Z\right\Vert ^{2}}=\frac{(k-1)!\omega_{k}(Z)}{\left\Vert Z\right\Vert ^{2k}%
.}.
\]
So the expression of $\left(  \partial\overline{\partial}g\right)  ^{d+k}$ is
reduced to%
\[
\left(  \partial\overline{\partial}g\right)  ^{d+k}=\frac{\left(  d+k\right)
!}{d!}\mu^{d}XY^{\prime}Y_{0}^{k-1}Y^{d}\left(  -\partial\overline{\partial
}\log N(z,z)\right)  ^{d}\wedge\frac{\omega_{k}(Z)}{\left\Vert Z\right\Vert
^{2k}.}.
\]
At $z=0$, we have $X=\left\Vert Z\right\Vert ^{2}$ and (\ref{EQN24}):
\[
\left.  \frac{\left(  -\partial\overline{\partial}\log N(z,z)\right)  ^{d}%
}{d!}\right\vert _{z=0}=\omega_{d}(z).
\]

Finally, we obtain%
\[
\frac{1}{(d+k)!}\left(  \partial\overline{\partial}g\right)  ^{d+k}%
(0,Z)=\mu^{d}X^{1-k}Y^{\prime}Y_{0}^{k-1}Y^{d}\omega,
\]
which is equivalent to (\ref{EQN45}).
\end{proof}

\begin{lemma}
\label{L9}Let $g$ be a $C^{2}$ function in $\widetilde{\Omega}$, which is a
solution of the Monge-Amp\`{e}re equation
\[
\det\left(  \frac{\partial^{2}g}{\partial z^{i}\partial\overline{z}^{j}%
}\right)  =\operatorname{e}^{(n+1)g}%
\]
and which generates an invariant form $\partial\overline{\partial}g$. Let
\begin{align*}
&  g(z,Z) =-\frac{\gamma+k\mu}{d+k+1}\log N(z,z)+h\left(  X\right)  ,\\
&  X =\frac{\left\Vert Z\right\Vert ^{2}}{N(z,z)^{\mu}}.
\end{align*}
Then $h$ satisfies on $]0,1[$ the differential equation
\begin{align*}
&  \left(  h^{\prime}(X)\right)  ^{k-1}\left(  \mu Xh^{\prime}(X)+\frac
{\gamma+k\mu}{d+k+1}\right)  ^{d}\left(  Xh^{\prime}(X)\right)  ^{\prime
}=\operatorname{e}^{(d+k+1)h},\\
&  Xh^{\prime}(X)\rightarrow0\qquad(X\rightarrow0).
\end{align*}
The boundary condition%
\[
g(z)\rightarrow\infty\qquad(z\rightarrow\partial\widetilde{\Omega})
\]
implies%
\[
h(X)\rightarrow\infty\qquad(X\rightarrow1).
\]

\end{lemma}

Let
\begin{align}
\beta &  =\frac{\gamma+k\mu}{\mu\left(  d+k+1\right)  },\label{EQN48}\\
Y  &  =Y_{0}+\beta=Xh^{\prime}(X)+\frac{\gamma+k\mu}{\mu\left(  d+k+1\right)
}. \label{EQN49}%
\end{align}
It results from the previous lemmas that if $g$ is a solution of the
Monge-Amp\`{e}re equation%
\[
\det\left(  \frac{\partial^{2}g}{\partial z^{i}\partial\overline{z}^{j}%
}\right)  =\operatorname{e}^{(d+k+1)g}%
\]
on $\widetilde{\Omega}$, with the boundary condition
\[
g(z)\rightarrow\infty\qquad(z\rightarrow\partial\widetilde{\Omega}),
\]
then the function $Y$ defined by (\ref{EQN49}) satisfies the differential
equation
\begin{equation}
\left(  \frac{Y-\beta}{X}\right)  ^{k-1}\left(  \mu Y\right)  ^{d}Y^{\prime
}=\operatorname{e}^{(d+k+1)h} \label{EQN60}%
\end{equation}
with the initial condition $Y(0)=\beta$.

Writing (\ref{EQN60}) as
\begin{equation}
\mu^{d}Y_{0}^{k-1}\left(  Y_{0}+\beta\right)  ^{d}Y_{0}^{\prime}%
=X^{k-1}\operatorname{e}^{(d+k+1)h} \label{EQN61}%
\end{equation}
shows that there exists a polynomial $R$ of degree $d+k$, with non negative
coefficients, such that the derivative of $R(Y_{0})$ is positive on $]0,1[$
and tends to $\infty$ when $X\rightarrow1$. So the function $Y_{0}$ is
strictly increasing and maps $[0,1[$ onto $[0,+\infty\lbrack$ and the function
$Y$ is strictly increasing and maps $[0,1[$ onto $[\beta,+\infty\lbrack$.

Taking logarithmic derivatives of both sides of (\ref{EQN61}), we get
\[
\frac{\left(  Y_{0}^{k-1}\left(  Y_{0}+\beta\right)  ^{d}Y_{0}^{\prime
}\right)  ^{\prime}}{Y_{0}^{k-1}\left(  Y_{0}+\beta\right)  ^{d}Y_{0}^{\prime
}}=(d+k+1)\frac{Y_{0}}{X}+\frac{k-1}{X},
\]
which is equivalent to
\[
X\left(  Y_{0}^{k-1}\left(  Y_{0}+\beta\right)  ^{d}Y_{0}^{\prime}\right)
^{\prime}=\left(  Y_{0}^{k-1}\left(  Y_{0}+\beta\right)  ^{d}Y_{0}^{\prime
}\right)  \left(  (d+k+1)Y_{0}+k-1\right)
\]
or to
\begin{equation}
\left(  XY_{0}^{k-1}\left(  Y_{0}+\beta\right)  ^{d}Y_{0}^{\prime}\right)
^{\prime}=\left(  (d+k+1)Y_{0}+k\right)  Y_{0}^{k-1}\left(  Y_{0}%
+\beta\right)  ^{d}Y_{0}^{\prime}. \label{EQN62}%
\end{equation}
Let $S$ be the polynomial of degree $d+1$ defined by%
\begin{equation}
T^{k}S(T)=\int_{0}^{T}\left(  (d+k+1)t+k\right)  t^{k-1}\left(  t+\beta
\right)  ^{d}\operatorname{d}t. \label{EQN62A}%
\end{equation}
Integrating (\ref{EQN62}) with the initial condition $Y_{0}(0)=0$ yields%
\[
XY_{0}^{k-1}\left(  Y_{0}+\beta\right)  ^{d}Y_{0}^{\prime}=Y_{0}^{k}S(Y_{0})
\]
or%
\[
X\left(  Y_{0}+\beta\right)  ^{d}Y_{0}^{\prime}=Y_{0}S(Y_{0}).
\]
Integrating by parts in (\ref{EQN62A}), we get%
\[
T^{k}S(T)=T^{k}(T+\beta)^{d+1}+k(1-\beta)\int_{0}^{T}t^{k}\left(
t+\beta\right)  ^{d}\operatorname{d}t;
\]
let
\begin{equation}
\int_{0}^{T}t^{k}\left(  t+\beta\right)  ^{d}\operatorname{d}t=T^{k+1}%
S_{1}(T). \label{EQN65}%
\end{equation}
Then $S_{1}$ is a polynomial of degree $d$, and
\[
(-1)^{k+d}\int_{0}^{-\beta}t^{k}\left(  t+\beta\right)  ^{d}\operatorname{d}%
t>0
\]
shows that $S_{1}(-\beta)\neq0$. Finally,
\begin{equation}
S(T)=(T+\beta)^{d+1}+k(1-\beta)TS_{1}(T),\quad S_{1}(-\beta)\neq0.
\label{EQN63}%
\end{equation}

\begin{lemma}
\label{L7}The differential equation
\begin{align*}
&  X\left(  Y_{0}+\beta\right)  ^{d}Y_{0}^{\prime} =Y_{0}S(Y_{0}),\\
&  Y_{0} \rightarrow\infty\qquad\qquad(X\rightarrow1),
\end{align*}
where the polynomial $S$ is defined by (\ref{EQN62A}), has a unique solution
\[
Y_{0}:[0,1[\rightarrow\lbrack0,+\infty\lbrack.
\]
This solution is $C^{\infty}$ at $0$.
\end{lemma}

\begin{theorem}
\label{T2}The generating function $g$ for the K\"{a}hler-Einstein metric of
\[
\widetilde{\Omega}=\widetilde{\Omega}_{k}(\mu)=\left\{  \left(  z,Z\right)
\in\Omega\times\mathbb{C}^{k}\mid\left\Vert Z\right\Vert ^{2}<N(z,z)^{\mu
}\right\}
\]
is given by
\[
g(z,Z)=-\frac{\gamma+k\mu}{d+k+1}\log N(z,z)+h\left(  \frac{\left\Vert
Z\right\Vert ^{2}}{N(z,z)^{\mu}}\right)  ,
\]
where
\[
\operatorname{e}^{(d+k+1)h}=\mu^{d}X^{1-k}Y_{0}^{\prime}Y_{0}^{k-1}\left(
Y_{0}+\beta\right)  ^{d}%
\]
and the function $Y_{0}:[0,1[\rightarrow\lbrack0,+\infty\lbrack$ satisfies%
\begin{align}
&  X\left(  Y_{0}+\beta\right)  ^{d}Y_{0}^{\prime}=Y_{0}S(Y_{0}),
\label{EQN64}\\
&  Y_{0}\rightarrow\infty\qquad\qquad(X\rightarrow1) \label{EQN64A}%
\end{align}
with
\[
T^{k}S(T)=\int_{0}^{T}\left(  (d+k+1)t+k\right)  t^{k-1}\left(  t+\beta
\right)  ^{d}\operatorname{d}t.
\]

\end{theorem}

In view of (\ref{EQN63}), the differential equation (\ref{EQN64}) may also be
written as
\begin{equation}
XY_{0}^{\prime}=Y_{0}(Y_{0}+\beta)+k(1-\beta)\frac{Y_{0}^{2}S_{1}(Y_{0}%
)}{\left(  Y_{0}+\beta\right)  ^{d}}, \label{EQN66}%
\end{equation}
where $S_{1}$ is the polynomial defined by (\ref{EQN65}).

\subsection{The critical exponent}

The expression (\ref{EQN63}) shows that the polynomial $S$ defined by
(\ref{EQN62A}) is divisible by $\left(  T+\beta\right)  ^{d}$ if and only if
$\beta=1$. In this case, $S=\left(  T+\beta\right)  ^{d+1}$. The equation
(\ref{EQN64}) is then%
\[
XY_{0}^{\prime}=Y_{0}\left(  Y_{0}+1\right)  .
\]
With the limit condition (\ref{EQN64A}), it integrates as
\[
-\log X=\int_{Y_{0}}^{\infty}\frac{1}{y(y+1)}\operatorname{d}y,
\]
which gives
\[
Y_{0}=\frac{X}{1-X}%
\]
and
\[
Y=Y_{0}+1=\frac{1}{1-X}.
\]
As
\[
\beta=\frac{\gamma+k\mu}{\mu\left(  d+k+1\right)  },
\]
the value of $\mu$ corresponding to $\beta=1$ is again%
\[
\mu_{0}=\frac{\gamma}{d+1},
\]
that is, the same value as for $k=1$.

For $\mu=\mu_{0}$, it is again possible to compute explicitly the
K\"{a}hler-Einstein metric of $\widetilde{\Omega}_{k}(\mu)$. Actually,%
\[
\operatorname{e}^{(d+k+1)h}=\mu^{d}X^{1-k}Y_{0}^{\prime}Y_{0}^{k-1}\left(
Y_{0}+\beta\right)  ^{d}%
\]
yields%
\begin{align*}
&  \operatorname{e}^{(d+k+1)h}=\mu_{0}^{d}\frac{1}{\left(  1-X\right)
^{d+k+1}},\\
&  h(X)=\frac{d}{d+k+1}\log\mu_{0}+\log\frac{1}{1-X}.
\end{align*}
As $\beta=1$, we have $\frac{\gamma+k\mu_{0}}{d+k+1}=\mu_{0},$
\begin{align*}
g(z,Z)  &  =-\mu_{0}\log N(z,z)+h\left(  \frac{\left\Vert Z\right\Vert ^{2}%
}{N(z,z)^{\mu_{0}}}\right) \\
&  =\frac{d}{d+k+1}\log\mu_{0}+\log\frac{1}{N(z,z)^{\mu_{0}}-\left\Vert
Z\right\Vert ^{2}}%
\end{align*}
and finally
\[
\partial\overline{\partial}g=-\partial\overline{\partial}\log\left(
N(z,z)^{\mu_{0}}-\left\Vert Z\right\Vert ^{2}\right)  .
\]

\section{The exceptional cases\label{Exceptional}}

Here we study the K\"{a}hler-Einstein metric and the Bergman metric of
$\widetilde{\Omega}_{k}\left(  \mu_{0}\right)  $, when $\Omega$ is one of the
two exceptional domains $\Omega_{V}$ and $\Omega_{VI}$ (see Appendix
\ref{ExceptionalDomains}).

\subsection{}

The \emph{Bergman kernel} of $\widetilde{\Omega}_{k}\left(
\mu\right)  $ has been computed (for all $\mu>0$) in
\cite{YinRoos2003}. If the polynomial $k\mapsto\chi(k\mu)$ (which
has degree $d=\dim\Omega$) is decomposed as
\begin{equation}
\frac{\chi(k\mu)}{\chi(0)}=\sum_{j=0}^{d}c_{\mu,j}\frac{\left(  k+1\right)
_{j}}{j!}, \label{VBComb}%
\end{equation}
where $\left(  k+1\right)  _{j}=\frac{\Gamma(k+j)}{\Gamma(k)}$ denotes the
raising factorial, let the function $F_{\chi,\mu}$ be defined by
\begin{equation}
F_{\chi,\mu}(t)=\sum_{j=0}^{d}c_{\mu,j}\left(  \frac{1}{1-t}\right)  ^{j}.
\label{VBRat}%
\end{equation}
Then the Bergman kernel $\widetilde{\mathcal{K}}_{k}(z,Z)$ of $\widetilde
{\Omega}_{k}\left(  \mu\right)  $ is given by the relations
\begin{align}
\widetilde{\mathcal{K}}_{k}(z,Z)  &  =\mathcal{L}_{k}\left(  z,\left\Vert
Z\right\Vert ^{2}\right)  ,\label{VBA}\\
\mathcal{L}_{k}(z,r)  &  =\frac{1}{k!}\frac{\partial^{k}}{\partial r^{k}%
}\mathcal{L}_{0}(z,r),\label{VBB}\\
\mathcal{L}_{0}\left(  z,r\right)   &  =\mathcal{K}(z)F_{\chi,\mu}\left(
\frac{r}{N(z,z)^{\mu}}\right)  , \label{VBC}%
\end{align}
where
\[
\mathcal{K}(z)=\frac{1}{\operatorname{vol}\Omega}\frac{1}{N\left(  z,z\right)
^{\gamma}}%
\]
is the Bergman kernel of $\Omega$.

Applying these results, we get
\begin{align*}
\mathcal{L}_{0}\left(  z,r\right)   &  =\frac{1}{\operatorname{vol}\Omega
}\frac{1}{N\left(  z,z\right)  ^{\gamma}}F_{\chi,\mu}\left(  \frac
{r}{N(z,z)^{\mu}}\right)  ,\\
\mathcal{L}_{k}(z,r)  &  =\frac{1}{k!}\frac{1}{\operatorname{vol}\Omega}%
\frac{1}{N\left(  z,z\right)  ^{\gamma+k\mu}}F_{\chi,\mu}^{(k)}\left(
\frac{r}{N(z,z)^{\mu}}\right)
\end{align*}
and%
\[
\widetilde{\mathcal{K}}_{k}(z,Z)=\frac{1}{k!}\frac{1}{\operatorname{vol}%
\Omega}\frac{1}{N\left(  z,z\right)  ^{\gamma+k\mu}}F_{\chi,\mu}^{(k)}(X),
\]
where
\[
X=\frac{\left\Vert Z\right\Vert ^{2}}{N(z,z)^{\mu}}.
\]
Hence the \emph{Bergman metric} of $\widetilde{\Omega}_{k}\left(  \mu\right)
$ is associated to the $(1,1)$ form%
\begin{equation}
\phi=-\left(  \gamma+k\mu\right)  \partial\overline{\partial}\log
N(z,z)+\partial\overline{\partial}\log F_{\chi,\mu}^{(k)}(X).. \label{B1}%
\end{equation}

If $\mu=\mu_{0}$, the critical exponent, then $\gamma+k\mu_{0}=\mu_{0}(d+k+1)$
and the Bergman metric of $\widetilde{\Omega}_{k}\left(  \mu_{0}\right)  $ is
associated to the $(1,1)$ form%
\begin{equation}
\phi_{0}=-\mu_{0}(d+k+1)\partial\overline{\partial}\log N(z,z)+\partial
\overline{\partial}\log F_{\chi,\mu_{0}}^{(k)}(X). \label{VBF}%
\end{equation}
For the critical exponent $\mu=\mu_{0}$, the K\"{a}hler-Einstein metric of
$\widetilde{\Omega}_{k}\left(  \mu_{0}\right)  $ is associated to the $(1,1)$
form
\begin{equation}
\Psi_{0}=\partial\overline{\partial}g=-\partial\overline{\partial}\log\left(
N(z,z)^{\mu_{0}}-\left\Vert Z\right\Vert ^{2}\right)  . \label{VBG}%
\end{equation}

\subsection{The exceptional case of dimension $16$\label{KernelExc16}}

See Appendix \ref{Exc16} for definitions, notations, and basic results.

Let
\[
V=\mathcal{M}_{2,1}(\mathbb{O}_{\mathbb{C}})=\left\{  \left(  a_{2}%
,a_{3}\right)  \mid a_{2},a_{3}\in\mathbb{O}_{\mathbb{C}}\right\}
\]
and consider the exceptional symmetric domain of dimension $16$%
\[
\Omega=\Omega_{V}=\left\{  x\in\mathcal{M}_{2,1}(\mathbb{O}_{\mathbb{C}}%
)\mid1-(x|x)+(x^{\sharp}\mid x^{\sharp})>0,\ 2-(x|x)>0\right\}  .
\]
Here $d=16$ and $\gamma=12$. The \emph{generic norm} is
\[
N(x,y)=1-(x|y)+(x^{\sharp}|y^{\sharp})\text{. }%
\]
The critical exponent is
\[
\mu_{0}=\frac{12}{17}.
\]
The inflated Hartogs domain is $\widetilde{\Omega}=\widetilde{\Omega}_{k}%
(\mu_{0})\subset\mathbb{C}^{16+k}$, consisting in $\left(  z,Z\right)
\in\mathbb{O}_{\mathbb{C}}\times\mathbb{C}^{k}$ which verify%
\begin{align*}
&  2-(z|z) >0,\\
&  1-(z|z)+(z^{\sharp} \mid z^{\sharp})>0,\\
&  \left\Vert Z\right\Vert ^{2} <\left(  1-(z|z)+(z^{\sharp}|z^{\sharp
})\right)  ^{12/17}.
\end{align*}
The \emph{K\"{a}hler-Einstein metric} of $\widetilde{\Omega}_{k}(\mu_{0})$ is
associated to the $(1,1)$-form%
\begin{equation}
\Psi_{0}=-\partial\overline{\partial}\log\left(  \left(  1-(z|z)+(z^{\sharp
}|z^{\sharp})\right)  ^{12/17}-\left\Vert Z\right\Vert ^{2}\right)  .
\label{VBH}%
\end{equation}

The \emph{Hua polynomial} $\chi$ is
\[
\chi(s)=(s+1)_{8}(s+4)_{8}.
\]
For $\mu=\mu_{0}=\frac{12}{17}$, the Bergman metric of $\widetilde{\Omega}%
_{k}\left(  \mu_{0}\right)  $ is associated to the $(1,1)$ form%
\begin{equation}
\phi_{0}=-\mu_{0}(d+k+1)\partial\overline{\partial}\log N(z,z)+\partial
\overline{\partial}\log F_{\chi,\mu_{0}}^{(k)}(X). \label{VBK}%
\end{equation}
Recall that
\[
F_{\chi,\mu_{0}}(t)=\sum_{j=0}^{d}c_{\mu_{0},j}\left(  \frac{1}{1-t}\right)
^{j},
\]
where the $c_{\mu,j}$ are determined by
\[
\frac{\chi(k\mu_{0})}{\chi(0)}=\sum_{j=0}^{16}c_{\mu_{0},j}\frac{\left(
k+1\right)  _{j}}{j!}.
\]

A computation with Maple gives
\[
\chi\left(  \frac{12}{17}s\right)  =\left(  \frac{12}{17}\right)  ^{16}%
\sum_{j=0}^{16}c_{j}(s+1)_{j},
\]
with%
\begin{align*}
&  c_{16}=1,\qquad c_{15}=0,\qquad c_{14}=\frac{595}{12},\qquad c_{13}%
=\frac{4165}{6},\qquad c_{12}=\frac{30042145}{3456},\\
&  c_{11}=\frac{14448385}{144},\qquad c_{10}=\frac{790269316375}%
{746496},\qquad c_{9}=\frac{1259425781075}{124416},\\
&  c_{8}=\frac{12447571001586875}{143327232},\qquad c_{7}=\frac
{2957566710311675}{4478976},\\
&  c_{6}=\frac{11300125622942496725}{2579890176},\qquad c_{5}=\frac
{10677213117341703625}{429981696},\\
&  c_{4}=\frac{65190770448545396318125}{557256278016},\qquad c_{3}%
=\frac{10209484788366056549125}{23219011584},\\
&  c_{2}=\frac{114818904611324955416375}{92876046336},\qquad c_{1}%
=\frac{35779252854815307462625}{15479341056},\\
&  c_{0}=\frac{33368892412222545303125}{15479341056}.
\end{align*}
This shows that \emph{all coefficients} $c_{\mu_{0},j}$ $(0\leq j\leq16)$ in
\[
F_{\chi,\mu_{0}}(t)=\sum_{j=0}^{16}c_{\mu_{0},j}\left(  \frac{1}{1-t}\right)
^{j}%
\]
\emph{are strictly positive, except }$c_{\mu_{0},15}=0$.

\subsection{The exceptional case of dimension $27$\label{KernelExc27}}

See Appendix \ref{Exc27} for definitions, notations, and basic results.

Let
\[
V=\mathcal{H}_{3}(\mathbb{O}_{\mathbb{C}}),
\]
the space of Cayley-Hermitian $3\times3$ matrices with entries in
$\mathbb{O}_{\mathbb{C}}$ and consider the exceptional symmetric domain
$\Omega=\Omega_{V}$ of dimension $27$, defined by the inequalities
\begin{gather*}
1-(z|z)+(z^{\sharp}|z^{\sharp})-\left\vert \det z\right\vert ^{2}>0,\\
3-2(z|z)+(z^{\sharp}|z^{\sharp})>0,\\
3-(z|z)>0.
\end{gather*}
The generic minimal polynomial is
\[
m(T,x,y)=T^{3}-(x|y)T^{2}+(x^{\sharp}|y^{\sharp})T-\det x\det\overline
{y}\text{.}%
\]
Here $d=27$ and $\gamma=18$. The \emph{generic norm }is
\[
N(x,y)=1-(x\mid y)+(x^{\#}\mid y^{\#})-\det x\det\overline{y}.
\]
The critical exponent is
\[
\mu_{0}=\frac{9}{14}..
\]
The inflated Hartogs domain is $\widetilde{\Omega}=\widetilde{\Omega}_{k}%
(\mu_{0})\subset\mathbb{C}^{16+k}$, consisting in $\left(  z,Z\right)
\in\mathbb{O}_{\mathbb{C}}\times\mathbb{C}^{k}$ which verify%
\begin{align*}
&  3-(z|z) >0,\\
&  3-2(z|z)+(z^{\sharp} \mid z^{\sharp})>0,\\
&  1-(z|z)+(z^{\sharp}|z^{\sharp})-\left\vert \det z\right\vert ^{2} >0,\\
&  \left\Vert Z\right\Vert ^{2} <\left(  1-(z|z)+(z^{\sharp}|z^{\sharp
})-\left\vert \det z\right\vert ^{2}\right)  ^{9/14}.
\end{align*}
The \emph{K\"{a}hler-Einstein metric} of $\widetilde{\Omega}_{k}(\mu_{0})$ is
associated to the $(1,1)$-form%
\begin{equation}
\Psi_{0}=-\partial\overline{\partial}\log\left(  \left(  1-(z|z)+(z^{\sharp
}|z^{\sharp})-\left\vert \det z\right\vert ^{2}\right)  ^{9/14}-\left\Vert
Z\right\Vert ^{2}\right)  .
\end{equation}

The \emph{Hua polynomial} $\chi$ is
\[
\chi(s)=(s+1)_{9}(s+5)_{9}(s+9)_{9}.
\]
For $\mu=\mu_{0}=\frac{9}{14}$, the Bergman metric of $\widetilde{\Omega}%
_{k}\left(  \mu_{0}\right)  $ is associated to the $(1,1)$ form%
\begin{equation}
\phi_{0}=-\mu_{0}(d+k+1)\partial\overline{\partial}\log N(z,z)+\partial
\overline{\partial}\log F_{\chi,\mu_{0}}^{(k)}(X),
\end{equation}
with
\[
F_{\chi,\mu_{0}}(t)=\sum_{j=0}^{27}c_{\mu_{0},j}\left(  \frac{1}{1-t}\right)
^{j},
\]
where the $c_{\mu,j}$ are determined by
\[
\frac{\chi(k\mu_{0})}{\chi(0)}=\sum_{j=0}^{27}c_{\mu_{0},j}\frac{\left(
k+1\right)  _{j}}{j!}.
\]

A computation with Maple gives
\[
\chi\left(  \frac{9}{14}s\right)  =\left(  \frac{9}{14}\right)  ^{27}%
\sum_{j=0}^{27}c_{j}(s+1)_{j},
\]
with%
\begin{align*}
c_{27}  &  =1,\qquad c_{26}=0,\qquad c_{25}=\frac{2275}{9},\\
c_{24}  &  =\frac{56875}{9},\qquad c_{23}=\frac{38591735}{243},\\
c_{21}  &  =15425515970150/3^{11},\\
c_{20}  &  =37061881356500/3^{9},\\
c_{19}  &  =184328710104188650/3^{14},\\
c_{18}  &  =3564334218619774600/3^{14},\\
c_{17}  &  =584735324681177419750/3^{16},\\
c_{16}  &  =10020732894163060819750/3^{16},\\
c_{15}  &  =352001611351295587864253500/3^{23},\\
c_{14}  &  =586664566244061492395923000/3^{21},\\
c_{13}  &  =1988637252859632373297511212000/3^{26},\\
c_{12}  &  =25672251717038124392289396301000/3^{26},\\
c_{11}  &  =8233663487061605972803486331644375/3^{29},\\
c_{10}  &  =89384793443821000370862374382625000/3^{29},\\
c_{9}  &  =1923754293102540042201539198326959366875/3^{36},\\
c_{8}  &  =209778908005712588859591649123533801875/3^{32},\\
c_{7}  &  =399192552377373476318550395682751432975625/3^{37},\\
c_{6}  &  =2728484170046421839052459199725228012518750/3^{37},\\
c_{5}  &  =47840351197962492631409316902739852226831250/3^{38},\\
c_{4}  &  =232468257762517753158460641861539626710125000/3^{38},\\
c_{3}  &  =73037107041363504672642146434776735686797778125/3^{42},\\
c_{2}  &  =23557400955895564936769134062033297681918662500/3^{40},\\
c_{1}  &  =409456797752799914624225389536137199476376953125/3^{42},\\
c_{0}  &  =394594700340674453245747775040231797415576953125/3^{42}.
\end{align*}
This shows again that \emph{all coefficients}\textit{ }$c_{\mu_{0},j}$ $(0\leq
j\leq27)$ in
\[
F_{\chi,\mu_{0}}(t)=\sum_{j=0}^{27}c_{\mu_{0},j}\left(  \frac{1}{1-t}\right)
^{j}%
\]
\emph{are strictly positive, except} $c_{\mu_{0},26}=0$.

\section{A conjecture\label{Conjecture}}

It has been shown above that for the critical exponent $\mu=\mu_{0}$, the
K\"{a}hler-Einstein metric of $\widetilde{\Omega}_{k}\left(  \mu_{0}\right)  $
is associated to the $(1,1)$ form
\[
\Psi_{0}=\partial\overline{\partial}g=-\partial\overline{\partial}\log\left(
N(z,z)^{\mu_{0}}-\left\Vert Z\right\Vert ^{2}\right)  .
\]
On the other hand, the Bergman metric of $\widetilde{\Omega}_{k}\left(
\mu\right)  $ is associated to the $(1,1)$ form%
\[
\phi_{0}=-\mu(d+k+1)\partial\overline{\partial}\log N(z,z)+\partial
\overline{\partial}\log F_{\chi,\mu}^{(k)}(X),
\]
with
\[
F_{\chi,\mu}(t)=\sum_{j=0}^{d}c_{\mu,j}\left(  \frac{1}{1-t}\right)  ^{j},
\]
where the $c_{\mu,j}$ are determined by
\begin{equation}
\frac{\chi(k\mu)}{\chi(0)}=\sum_{j=0}^{d}c_{\mu,j}\frac{\left(  k+1\right)
_{j}}{j!}. \label{VBComb1}%
\end{equation}

\begin{conjecture}
Let $\Omega$ be an irreducible circled bounded symmetric domain of genus
$\gamma$ and dimension $d$, $\mu_{0}=\frac{\gamma}{d+1}$ its critical
exponent, $\chi$ its Hua polynomial. The coefficients $c_{\mu,j}$ in
(\ref{VBComb1}) are all strictly positive if and only if%
\[
\mu<\mu_{0}.
\]
For $\mu=\mu_{0}$, all coefficients $c_{\mu,j}$ in (\ref{VBComb1}) are strictly
positive, except $c_{\mu,d-1}=0$ and except for the rank $1$ type $I_{1,n}$
(where $c_{\mu,d}>0$ and $c_{\mu,j}=0$ for all $j<d=n$).
\end{conjecture}

The values of the critical exponent are

\medskip%

\begin{tabular}
[c]{lllllll}%
Type & $I_{m,n}$ & $II_{n}$ & $III_{n}$ & $IV_{n}$ & $V$ & $VI$\\
$\mu_{0}$ & $\frac{m+n}{mn+1}$ & $\frac{2}{n+\frac{2}{n-1}}$ &
$\frac
{2}{n+\frac{2}{n+1}}$ & $\frac{n}{n+1}$ & $\frac{12}{17}$ & $\frac{9}{14}$%
\end{tabular}

\medskip

\noindent We have always $\mu_{0}<1$, except in the rank $1$ case $I_{1,n}$.
For $\mu=1$, the signs of $c_{\mu,j}$ are alternating, starting with
$c_{\mu,d}>0$ and ending with $c_{\mu,j}=0$ for $j<j_{0}$, where $j_{0}$ is a
positive integer depending on $\Omega$.

\begin{remark}
The conjecture has been checked with help of computer algebra software in many
significant cases:

\begin{itemize}
\item for $\mu=\mu_{0}$ and the types $I_{3,3}$, $IV_{3}$, $IV_{4}$, $IV_{6}$,
$V$, $VI$;

\item for type $V$ and all values of $\mu$.
\end{itemize}
\end{remark}

\begin{remark}
As the function $F_{\chi,\mu}$ is related to the Bergman kernel of
$\widetilde{\Omega}_{k}\left(  \mu\right)  $, all derivatives $F_{\chi,\mu
}^{(k)}$, $k>0$, of this function are strictly positive on $[0,1[$ for all
$\mu>0$.
\end{remark}

\begin{remark}
If the conjecture is true, it would help to compare the Bergman metric and the
K\"{a}hler-Einstein metric of $\widetilde{\Omega}_{k}\left(  \mu\right)  $ for
$\mu=\mu_{0}$. The exponent $\mu_{0}$ seems also to be a limit case for other
comparisons: in \cite{YinWangZhao2001}, a comparison theorem is given between
the K\"{a}hler-Einstein metric and the Kobayashi metric of $\widetilde{\Omega
}_{1}\left(  \mu\right)  $ when $\Omega$ is a symmetric domain of type
$I_{m,n}$ and $\mu<\mu _{0}$.
\end{remark}

\appendix

\section{Bounded symmetric domains and Jordan triple systems}

Hereunder we give a review of properties of the Jordan triple structure
associated to a complex bounded symmetric domain (see \cite{Loos1977},
\cite{Roos1999}).

\subsection{Jordan triple system associated to a bounded symmetric
domain\label{JTS}}

Let $\Omega$ be an irreducible bounded circled homogeneous domain in a complex
vector space $V$. Let $K$ be the identity component of the (compact) Lie group
of (linear) automorphisms of $\Omega$ leaving $0$ fixed. Let $\omega$ be a
volume form on $V$, invariant by $K$ and by translations. Let $\mathcal{K}$ be
the Bergman kernel of $\Omega$ with respect to $\omega$, that is, the
reproducing kernel of the Hilbert space $H^{2}(\Omega,\omega)=\mathrm{Hol}%
(\Omega)\cap L^{2}(\Omega,\omega)$. The Bergman metric at $z\in\Omega$ is
defined by
\[
h_{z}(u,v)=\partial_{u}\overline{\partial}_{v}\log\mathcal{K}(z).
\]
The \emph{Jordan triple product} on $V$ is characterized by
\[
h_{0}(\{uvw\},t)=\partial_{u}\overline{\partial}_{v}\partial_{w}%
\overline{\partial}_{t}\log\mathcal{K}(z)\left\vert _{z=0}\right.  .
\]
The triple product $(x,y,z)\mapsto\{xyz\}$ is complex bilinear and symmetric
with respect to $(x,z)$, complex antilinear with respect to $y$. It satisfies
the\emph{\ Jordan identity}
\[
\{xy\{uvw\}\}-\{uv\{xyw\}\}=\{\{xyu\}vw\}-\{u\{vxy\}w\}.
\]
The space $V$ endowed with the triple product $\{xyz\}$ is called a
\emph{(Hermitian) Jordan triple system}. For $x,y,z\in V$, denote by $D(x,y)$
and $Q(x,z)$ the operators defined by
\[
\{xyz\}=D(x,y)z=Q(x,z)y.
\]
The Bergman metric at $0$ is related to $D$ by
\[
h_{0}(u,v)=\operatorname*{tr}D(u,v).
\]
A Jordan triple system is called \emph{Hermitian positive }if
$(u|v)=\operatorname*{tr}D(u,v)$ is positive definite. As the Bergman metric
of a bounded domain is always definite positive, the Jordan triple system
associated to a bounded symmetric domain is Hermitian positive.

The \emph{quadratic representation }
\[
Q:V\longrightarrow\operatorname*{End}{}_{\mathbb{R}}(V)
\]
is defined by $Q(x)y=\frac{1}{2}\{xyx\}$. The following fundamental identity
for the quadratic representation is a consequence of the Jordan identity:
\[
Q(Q(x)y)=Q(x)Q(y)Q(x).
\]
The \emph{Bergman operator} $B$ is defined by
\[
B(x,y)=I-D(x,y)+Q(x)Q(y),
\]
where $I$ denotes the identity operator in $V$. It is also a consequence of
the Jordan identity that the following fundamental identity holds for the
Bergman operator:
\[
Q(B(x,y)z)=B(x,y)Q(z)B(y,x).
\]
The Bergman operator gets its name from the following property:
\[
h_{z}\left(  B(z,z)u,v\right)  =h_{0}(u,v)\quad(z\in\Omega;\ u,v\in V).
\]
If $\Phi\in\left(  \operatorname*{Aut}\Omega\right)  _{0}$, the identity
component of the automorphism group of $\Omega$, the relation
\begin{equation}
B(\Phi x,\Phi y)=\operatorname{d}\Phi(x)\circ B(x,y)\circ\operatorname{d}%
\Phi(y)^{\ast} \label{K12}%
\end{equation}
holds for $x,y\in\Omega$, where $~^{\ast}$ denotes the adjoint with respect to
the Hermitian metric $h_{0}$. As a consequence, the Bergman kernel of $\Omega$
is given by
\begin{equation}
\mathcal{K}(z)=\frac{1}{\operatorname*{vol}\Omega}\frac{1}{\det B(z,z)}.
\label{K1}%
\end{equation}

The \emph{quasi-inverse }$x^{y}$ is defined, for each pair $\left(
x,y\right)  $ such that $B(x,y)$ is invertible, by
\[
x^{y}=B\left(  x,y\right)  ^{-1}\left(  x-Q(x)y\right)  .
\]

\subsection{Spectral theory}

An Hermitian positive Jordan triple system is always \emph{semi-simple}, that
is, the direct sum of a finite family of simple subsystems, with
compo\-nent-wise triple product.

As the domain $\Omega$ is assumed to be irreducible, the associated Jordan
triple system $V$ is \emph{simple}, that is $V$ is not the direct sum of two
non trivial subsystems.

An \emph{automorphism} $f:V\rightarrow V$ of the Jordan triple system $V$ is a
complex linear isomorphism preserving the triple product~:
$f\{u,v,w\}=\{fu,fv,fw\}$. The automorphisms of $V$ form a group, denoted
$\operatorname{Aut}V$, which is a compact Lie group; we will denote by $K$ its
identity component.

An element $c\in V$ is called \emph{tripotent} if $\{ccc\}=2c$. If $c$ is a
tripotent, the operator $D(c,c)$ annihilates the polynomial $T(T-1)(T-2)$.

Let $c$ be a tripotent. The decomposition $V=V_{0}(c)\oplus V_{1}(c)\oplus
V_{2}(c) $, where $V_{j}(c)$ is the eigenspace $V_{j}(c)=\left\{  x\in
V\;;\;D(c,c)x=jx\right\}  $, is called the \emph{Peirce decomposition} of $V$
(with respect to the tripotent $c$).

Two tripotents $c_{1}$ and $c_{2}$ are called \emph{orthogonal} if
$D(c_{1},c_{2})=0$. If $c_{1}$ and $c_{2}$ are orthogonal tripotents, then
$D(c_{1},c_{1)}$ and $D(c_{2},c_{2})$ commute and $c_{1}+c_{2}$ is also a tripotent.

A non zero tripotent $c$ is called \emph{primitive} if it is not the sum of
non zero orthogonal tripotents. A tripotent $c$ is \emph{maximal} if there is
no non zero tripotent orthogonal to $c$. The set of maximal tripotents is
equal to the \emph{Shilov boundary} of the domain $\Omega$.

A \emph{frame} of $V$ is a maximal sequence $(c_{1},\ldots,c_{r})$ of pairwise
orthogonal primitive tripotents. The frames of $V$ form a manifold
$\mathcal{F}$, which is called the \emph{Satake-Furstenberg boundary} of
$\Omega$.

Let $\mathbf{c}=(c_{1},\ldots,c_{r})$ be a frame. For $0\leq i\leq j\leq r$,
let
\[
V_{ij}(\mathbf{c})=\left\{  x\in V\mid D(c_{k},c_{k})x=(\delta_{i}^{k}%
+\delta_{j}^{k})x,\;1\leq k\leq r\right\}  \text{~: }%
\]
the decomposition $V=\bigoplus_{0\leq i\leq j\leq r}V_{ij}(\mathbf{c})$ is
called the \emph{simultaneous Peirce decomposition} with respect to the frame
$\mathbf{c}$.

Let $V$ be a simple Hermitian positive Jordan triple system. Then there exist
frames for $V$. All frames have the same number of elements, which is the
\emph{rank} $r$ of $V$. The subspaces $V_{ij}=V_{ij}(\mathbf{c})$ of the
simultaneous Peirce decomposition have the following properties: $V_{00}=0$~;
$V_{ii}=\mathbb{C}e_{i}$ ($0<i$); all $V_{ij}$'s ($0<i<j$) have the same
dimension $a$; all $V_{0i}$'s ($0<i$) have the same dimension $b$.

The \emph{numerical invariants} of $V$ (or of $\Omega$) are the rank $r$ and
the two integers
\begin{gather*}
a=\dim V_{ij}\;\;(0<i<j)\text{,}\\
b=\dim V_{0i}\;\;(0<i)\text{.}%
\end{gather*}
The \emph{genus} of $V$ is the number $\gamma$ defined by
\[
\gamma=2+a(r-1)+b\text{;}%
\]
the genus is generally denoted by $g$ or $p$, but we denote it here by
$\gamma$ as $g$ stands for the generating function of the K\"{a}hler-Einstein
metric. The HPJTS $V$ and the domain $\Omega$ are said to be of \emph{tube
type} if $b=0$.

Let $V$ be a simple Hermitian positive Jordan triple system. Then any $x\in V$
can be written in a unique way
\begin{equation}
x=\lambda_{1}c_{1}+\lambda_{2}c_{2}+\cdots+\lambda_{p}c_{p}\text{,} \label{A1}%
\end{equation}
where $\lambda_{1}>\lambda_{2}>\cdots>\lambda_{p}>0$ and $c_{1},c_{2}%
\ldots,c_{p}$ are pairwise orthogonal tripotents. The element $x$ is
\emph{regular} iff $p=r$; then $(c_{1},c_{2},\ldots,c_{r})$ is a frame of $V$.
The decomposition (\ref{A1}) is called the \emph{spectral decomposition} of
$x$.

\subsection{The generic minimal polynomial\label{GenericNorm}}

Let $V$ be a Jordan triple system of rank $r$. There exist polynomials
$m_{1},\ldots,m_{r}$ on $V\times\overline{V}$, homogeneous of respective
bidegrees $(1,1),\ldots,(r,r)$, such that for each regular $x\in V$, the
polynomial
\[
m(T,x,y)=T^{r}-m_{1}(x,y)T^{r-1}+\cdots+(-1)^{r}m_{r}(x,y)
\]
satisfies
\[
m(T,x,x)=\prod_{i=1}^{r}(T-\lambda_{i}^{2})\text{, }%
\]
where $x=\lambda_{1}c_{1}+\lambda_{2}c_{2}+\cdots+\lambda_{r}c_{r}$ is the
spectral decomposition of $x$. Here $\overline{V}$ denotes the space $V$ with
the conjugate complex structure. The polynomial
\[
m(T,x,y)=T^{r}-m_{1}(x,y)T^{r-1}+\cdots+(-1)^{r}m_{r}(x,y)
\]
is called the \emph{generic minimal polynomial} of $V$ (at $(x,y)$). The
(inhomogeneous) polynomial $N:V\times\overline{V}\rightarrow\mathbb{C}$
defined by
\[
N(x,y)=m(1,x,y)
\]
is called the \emph{generic norm}. The following identities hold:
\begin{align}
\det B(x,y)  &  =N(x,y)^{\gamma},\label{K2}\\
\operatorname*{tr}D(x,y)  &  =\gamma\,m_{1}(x,y)\text{.}\nonumber
\end{align}

\subsection{The spectral norm}

Let $V$ be an HPJTS. The map $x\mapsto\lambda_{1}$, where $x=\lambda_{1}%
c_{1}+\lambda_{2}c_{2}+\cdots+\lambda_{p}c_{p}$ is the spectral decomposition
of $x$ ($\lambda_{1}>\lambda_{2}>\cdots>\lambda_{p}>0)$ is a norm on $V$,
called the \emph{spectral norm}. The \emph{bounded symmetric domain} $\Omega$
is the unit ball of $V$ for the \emph{spectral norm}. It is also characterized
by the set of polynomial inequalities
\[
\left.  \frac{\partial^{j}}{\partial T^{j}}m(T,x,x)\right\vert _{T=1}%
>0,\qquad0\leq j\leq r-1.
\]

\subsection{The Hua integral}

Let $\Omega$ be an irreducible bounded circled homogeneous domain, $N$ its
generic norm. The \emph{Hua integral} relative to $\Omega$ is
\[
\int_{\Omega}N(x,x)^{s}\alpha^{n},
\]
where $\alpha=\frac{\operatorname{i}}{2\pi}\partial\overline{\partial}m_{1}$.
It converges for $s>-1$ and has be computed for the four classical series of
symmetric domains by Hua Lookeng in \cite{Hua1963}.

For a general irreducible bounded circled homogeneous domain with numerical
invariants $a$, $b$ and rank $r$, the Hua integral is given by
\[
\int_{\Omega}N(x,x)^{s}\alpha^{n}=\frac{\chi(0)}{\chi(s)}\int_{\Omega}%
\alpha^{n},
\]
where $\chi$ is the \emph{Hua polynomial}
\[
\chi(s)=\prod_{j=1}^{r}\left(  s+1+(j-1)\frac{a}{2}\right)  _{1+b+\left(
r-j\right)  a}%
\]
(see \cite{YinRoos2003}, Theorem 2.3). See the tables below for the expression
of this polynomial for the six cases of bounded irreducible symmetric domains.

\section{Tables for bounded symmetric domains\label{Tables}}

The following examples exhaust the list of simple Hermitian positive Jordan
triple systems (see \cite{Loos1975}). The HPJTS occurring in the four infinite
series $I_{p,q}$, $II_{n}$, $III_{n}$, $IV_{n}$ are called \emph{classical};
the two HPJTS of type $V$ and $VI$ are called \emph{exceptional}. There is
some overlapping between the classical series, due to a finite number of
isomorphisms in low dimension. We give hereunder for each type:

\begin{itemize}
\item the definition of the space $V$, its Jordan triple product, and the
corresponding bounded circled homogeneous domain;

\item the generic norm;

\item the numerical invariants $r,a,b,\gamma=2+a(r-1)+b$;

\item the Hua polynomial
\[
\chi(s)=\prod_{j=1}^{r}\left(  s+1+(j-1)\frac{a}{2}\right)  _{1+b+\left(
r-j\right)  a}.
\]

\end{itemize}

\subsection{Type I$_{m,n}$ $\left(  1\leq m\leq n\right)  $}

$V=\mathcal{M}_{m,n}(\mathbb{C})$ (space of $m\times n$ matrices with complex
entries), endowed with the triple product
\[
\{xyz\}=x^{t}\overline{y}z+z^{t}\overline{y}x.
\]
The domain $\Omega$ is the set of $m\times n$ matrices $x$ such that
$I_{m}-x^{t}\overline{x}$ is definite positive. The generic minimal polynomial
is
\[
m(T,x,y)=\operatorname*{Det}(TI_{m}-x^{t}\overline{y}),
\]
where $\operatorname*{Det}$ is the usual determinant of square matrices. The
numerical invariants are $r=m$, $a=2$, $b=n-m$, $\gamma=m+n$. These HPJTS are
of tube type only for $m=n$.

The polynomial $\chi$ is%
\[
\chi(s)=\prod_{j=1}^{m}\left(  s+j\right)  _{m+n+1-2j}=\prod_{j=1}%
^{m}(s+j)_{n}.
\]

\subsection{Type II$_{n}$ $\left(  n\geq2\right)  $}

$V=\mathcal{A}_{n}(\mathbb{C})$ (space of $n\times n$ alternating matrices)
with the same triple product as for Type I. The domain $\Omega$ is the set of
$n\times n$ alternating matrices $x$ such that $I_{n}+x\overline{x}$ is
definite positive.

\subsubsection{Type II$_{2p}$ ($n=2p$ even)}

The generic minimal polynomial is here given by
\[
m(T,x,y)^{2}=\operatorname*{Det}(TI_{n}+x\overline{y}).
\]
The numerical invariants are $r=\frac{n}{2}=p$, $a=4$, $b=0$, $\gamma=2\left(
n-1\right)  $; these HPJTS are of tube type.

The polynomial $\chi$ is
\[
\chi(s)=\prod_{j=1}^{p}\left(  s+2j-1\right)  _{1+4\left(  p-j\right)  }%
=\prod_{j=1}^{p}(s+2j-1)_{2p-1}.
\]

\subsubsection{Type II$_{2p+1}$ ($n=2p+1$ odd)}

The generic minimal polynomial is given by
\[
Tm(T,x,y)^{2}=\operatorname*{Det}(TI_{n}+x\overline{y}).
\]
The numerical invariants are $r=\left[  \frac{n}{2}\right]  =p$, $a=4$, $b=2$,
$\gamma=2(n-1)$; these HPJTS are not of tube type.

The polynomial $\chi$ is%
\[
\chi(s)=\prod_{j=1}^{p}\left(  s+2j-1\right)  _{3+4\left(  p-j\right)  }%
=\prod_{j=1}^{p}(s+2j-1)_{2p+1}.
\]

\subsection{Type III$_{n}$ $\left(  n\geq1\right)  $}

$V=\mathcal{S}_{n}(\mathbb{C})$ (space of $n\times n$ symmetric matrices) with
the same triple product as for Type I. The domain $\Omega$ is the set of
$n\times n$ symmetric matrices $x$ such that $I_{n}-x\overline{x}$ is definite
positive. The generic minimal polynomial is
\[
m(T,x,y)=\operatorname*{Det}(TI_{n}-x\overline{y}).
\]
The numerical invariants are $r=n$, $a=1$, $b=0$, $\gamma=n+1$. These HPJTS
are of tube type.

The polynomial $\chi$ is
\[
\chi(s)=\prod_{j=1}^{n}\left(  s+\frac{j+1}{2}\right)  _{1+n-j}.
\]

\subsection{Type IV$_{n}$ $\left(  n\neq2\right)  $}

$V=\mathbb{C}^{n}$ with the quadratic operator defined by
\[
Q(x)y=q(x,\overline{y})x-q(x)\overline{y},
\]
where $q(x)=\sum x_{i}^{2},$ $q(x,y)=2\sum x_{i}y_{i}$. The domain $\Omega$ is
the set of points $x\in\mathbb{C}^{n}$ such that
\[
1-q(x,\overline{x})+\left\vert q(x)\right\vert ^{2}>0,\qquad2-q(x,\overline
{x})>0.
\]
The generic minimal polynomial is
\[
m(T,x,y)=T^{2}-q(x,\overline{y})+q(x)q(\overline{y}).
\]
The numerical invariants are $r=2$, $a=n-2$, $b=0$, $\gamma=n$. These HPJTS
are of tube type.

The polynomial $\chi$ is
\[
\chi(s)=\left(  s+1\right)  _{n-1}\left(  s+\frac{n}{2}\right)  .
\]

\subsection{Type V}

$V=\mathcal{M}_{2,1}(\mathbb{O}_{\mathbb{C}})$, the subspace of $\mathcal{H}%
_{3}(\mathbb{O}_{\mathbb{C}})$ consisting in matrices of the form
\[
\left(
\begin{array}
[c]{ccc}%
0 & a_{3} & \tilde{a}_{2}\\
\tilde{a}_{3} & 0 & 0\\
a_{2} & 0 & 0
\end{array}
\right)
\]
with the same quadratic operator as for type VI (see below). Here $\tilde{a}$
denotes the Cayley conjugate of $a\in\mathbb{O}_{\mathbb{C}}$. The generic
minimal polynomial is
\[
m(T,x,y)=T^{2}-(x|y)T+(x^{\sharp}|y^{\sharp})\text{. }%
\]
The domain $\Omega$ is the \textquotedblleft exceptional domain of dimension
$16$\textquotedblright\ defined by
\[
1-(x|x)+(x^{\sharp}\mid x^{\sharp})>0,\qquad2-(x|x)>0.
\]
The numerical invariants are $r=2$, $a=6$, $b=4$, $\gamma=12$. This HPJTS is
not of tube type.

The polynomial $\chi$ is
\[
\chi(s)=\left(  s+1\right)  _{11}\left(  s+4\right)  _{5};
\]
it can also be written
\[
\chi(s)=(s+1)_{8}(s+4)_{8}.
\]

\subsection{Type VI}

$V=\mathcal{H}_{3}(\mathbb{O}_{\mathbb{C}})$, the space of $3\times3$ matrices
with entries in the space $\mathbb{O}_{\mathbb{C}}$ of octonions over
$\mathbb{C}$, which are Hermitian with respect to the Cayley conjugation; the
quadratic operator is defined by
\[
Q(x)y=(x|y)x-x^{\sharp}\times\overline{y},
\]
where $\times$ denotes the Freudenthal product, $x^{\sharp}$ the adjoint
matrix in $\mathcal{H}_{3}(\mathbb{O}_{\mathbb{C}})$ and $(x|y)$ the standard
Hermitian product in $\mathcal{H}_{3}(\mathbb{O}_{\mathbb{C}})$ (see Appendix
\ref{Exc27}). The domain $\Omega$ is the \textquotedblleft exceptional domain
of dimension $27$\textquotedblright\ defined by
\begin{gather*}
1-(x|x)+(x^{\sharp}|x^{\sharp})-\left\vert \det x\right\vert ^{2}>0,\\
3-2(x|x)+(x^{\sharp}|x^{\sharp})>0,\\
3-(x|x)>0.
\end{gather*}
The generic minimal polynomial is
\[
m(T,x,y)=T^{3}-(x|y)T^{2}+(x^{\sharp}|y^{\sharp})T-\det x\det\overline
{y}\text{, }%
\]
where $\det$ denotes the determinant in $\mathcal{H}_{3}(\mathbb{O}%
_{\mathbb{C}})$. The numerical invariants are $r=3$, $a=8$, $b=0$, $\gamma
=18$. This HPJTS is of tube type.

The polynomial $\chi$ is
\begin{align*}
\chi(s)  &  =(s+1)_{17}(s+5)_{9}(s+9)\\
&  =(s+1)_{9}(s+5)_{9}(s+9)_{9}.
\end{align*}

\section{The exceptional bounded symmetric domains\label{ExceptionalDomains}}

In this appendix, we recall without proofs the construction and
the main results about the two exceptional bounded symmetric
domains. For details and proofs, see \cite{Roos1992}.

\subsection{The exceptional Jordan triple system $H_{3}(\mathbb{O}%
_{\mathbb{C}})$ and the exceptional bounded symmetric domain of dimension
$27$\label{Exc27}}

Let $\mathbb{O}_{\mathbb{C}}$ denote the $8$-dimen\-sional algebra of complex
octonions, with Cayley conjugation $a\mapsto\tilde{a}$, trace $t(a)=a+\tilde
{a}$, Cayley norm $n(a)=a\tilde{a}$. If $\mathbb{O}$ denotes the real Cayley
division algebra, we consider $\mathbb{O}_{\mathbb{C}}$ as its
complexification $\mathbb{O}_{\mathbb{C}}=\mathbb{C\otimes}_{\mathbb{R}%
}\mathbb{O}$; the complex conjugate of $a=\lambda\alpha$ ($\lambda
\in\mathbb{C}$, $\alpha\in\mathbb{O}$) is $\overline{a}=\overline{\lambda
}\alpha$. (Recall that the algebras $\mathbb{O}$ and $\mathbb{O}_{\mathbb{C}}$
are neither commutative nor associative). The (complex bilinear) scalar
product of $a,b\in\mathbb{O}_{\mathbb{C}}$ is defined by
\begin{equation}
\left(  a:b\right)  =a\widetilde{b}+\widetilde{a}b=\widetilde{b}%
a+b\widetilde{a}; \label{eq2/1b}%
\end{equation}
the Hermitian scalar product of $a,b\in\mathbb{O}_{\mathbb{C}}$ is defined by
\begin{equation}
\left(  a\mid b\right)  =\left(  a:\overline{b}\right)  =a\widetilde
{\overline{b}}+\widetilde{a}\overline{b}=\widetilde{\overline{b}}%
a+\overline{b}\widetilde{a}. \label{eq2/1c}%
\end{equation}

We denote by $H_{3}(\mathbb{O}_{\mathbb{C}})$ the $\mathbb{C}$-vector space
(with the natural operations) of $3\times3$ matrices with entries in
$\mathbb{O}_{\mathbb{C}}$, which are \emph{Hermitian with respect to the
Cayley conjugation in }$O_{\mathbb{C}}$\emph{. }An element $a\in
H_{3}(\mathbb{O}_{\mathbb{C}})$ will be written
\begin{equation}
a=\left(
\begin{array}
[c]{lll}%
\alpha_{1} & a_{3} & \widetilde{a_{2}}\\
\widetilde{a_{3}} & \alpha_{2} & a_{1}\\
a_{2} & \widetilde{a_{1}} & \alpha_{3}%
\end{array}
\right)  , \label{eq2/1}%
\end{equation}
with $\alpha_{1},\alpha_{2},\alpha_{3}\in\mathbb{C}$ and $a_{1},a_{2},a_{3}%
\in\mathbb{O}_{\mathbb{C}}$. Instead of (\ref{eq2/1}), we will also write
\begin{equation}
a=\sum_{j=1}^{3}\alpha_{j}e_{j}+\sum_{j=1}^{3}F_{j}(a_{j}), \label{eq2/1A}%
\end{equation}
with the obvious definitions for $e_{j}$ and $F_{j}(a_{j})$. The vector space
$H_{3}(\mathbb{O}_{\mathbb{C}})$ decomposes into the direct sum
\begin{equation}
H_{3}(\mathbb{O}_{C})=\mathbb{C}e_{1}\oplus\mathbb{C}e_{2}\oplus
\mathbb{C}e_{3}\oplus\mathcal{F}_{1}\oplus\mathcal{F}_{2}\oplus\mathcal{F}%
_{3}, \label{eq2/2}%
\end{equation}
where $\mathcal{F}_{j}=\left\{  F_{j}(a);a\in\mathbb{O}_{\mathbb{C}}\right\}
$. The subspaces $\mathcal{F}_{j}$ are $8$-dimensional and
\[
\dim_{\mathbb{C}}H_{3}(\mathbb{O}_{C})=27.
\]
On $H_{3}(\mathbb{O}_{\mathbb{C}})$, define a bilinear form by
\begin{equation}
(a:b)=\sum_{j=1}^{3}\alpha_{j}\beta_{j}+\sum_{j=1}^{3}(a_{j}:b_{j})
\label{eq2/3}%
\end{equation}
for $a=\sum_{j=1}^{3}\alpha_{j}e_{j}+\sum_{j=1}^{3}F_{j}(a_{j})$,
$b=\sum_{j=1}^{3}\beta_{j}e_{j}+\sum_{j=1}^{3}F_{j}(b_{j})$. Here
$(a_{j}:b_{j})$ denotes the scalar product (\ref{eq2/17}) in $\mathbb{O}%
_{\mathbb{C}}$. The form defined by (\ref{eq2/3}) is clearly non-singular and
the decomposition (\ref{eq2/2}) is orthogonal with respect to it. We will
refer to $(a:b)$ as the (complex bilinear) \emph{scalar product} of $a$ and
$b$ in $H_{3}(\mathbb{O}_{\mathbb{C}})$.

The \emph{adjoint} $a^{\#}$ of an element
\[
a=\sum_{j=1}^{3}\alpha_{j}e_{j}+\sum_{j=1}^{3}F_{j}(a_{j})\in H_{3}%
(\mathbb{O}_{\mathbb{C}})
\]
is defined by
\begin{equation}
a^{\#}=\sum_{i}\left(  \alpha_{j}\alpha_{k}-n(a_{i})\right)  e_{i}+\sum
_{i}\widetilde{F_{i}}\left(  a_{j}a_{k}-\alpha_{i}\widetilde{a_{i}}\right)  .
\label{eq2/4}%
\end{equation}
In (\ref{eq2/4}) and below, $\sum_{i}$ means $\sum_{i=1}^{3}$ and $j,k$ are
defined by $(i,j,k)$ being an \emph{even} permutation of $(1,2,3)$;
$\widetilde{F_{i}}(c)$ stands for $F_{i}(\widetilde{c})$. The symmetric
bilinear map, associated to the quadratic map $a\mapsto a^{\#}$ by%
\[
a\times b=(a+b)^{\#}-a^{\#}-b^{\#},\quad a\times a=2a^{\#},
\]
is called the \emph{Freudenthal product} in $H_{3}(\mathbb{O}_{\mathbb{C}})$.
It follows directly from the definitions that
\begin{equation}
a\times b=\sum_{i}\left(  \alpha_{j}\beta_{k}+\alpha_{k}\beta_{j}-(a_{i}%
:b_{i})\right)  e_{i}+\sum_{i}\widetilde{F_{i}}\left(  a_{j}b_{k}+b_{j}%
a_{k}-\alpha_{i}\widetilde{b_{i}}-\beta_{i}\widetilde{a_{i}}\right)  .
\label{eq2/26}%
\end{equation}
In particular, we have%
\begin{gather}
e_{i}\times e_{i}=0,\quad e_{i}\times e_{j}=e_{k},\nonumber\\
e_{i}\times F_{i}(b)=-F_{i}(b),\quad e_{i}\times F_{j}(b)=0,\quad
\label{eq2/27}\\
F_{i}(a)\times F_{i}(b)=-(a:b)e_{i},\quad F_{i}(a)\times F_{j}(b)=\widetilde
{F_{k}}(ab).\nonumber
\end{gather}
In these relations, $(i,j,k)$ is always an even permutation of $(1,2,3)$. The
Freudenthal product verifies
\begin{equation}
\left(  a\times b:c\right)  =\left(  a:b\times c\right)  \qquad(a,b,c\in
H_{3}(\mathbb{O}_{\mathbb{C}})). \label{eq2/5}%
\end{equation}

Let $T$ denote the trilinear symmetric form on $H_{3}(\mathbb{O}_{\mathbb{C}%
})$ defined by
\[
T(a,b,c)=\left(  a\times b:c\right)  .
\]
The \emph{determinant} in $H_{3}(\mathbb{O}_{\mathbb{C}})$ is the associated
polynomial of degree $3$, defined by
\begin{equation}
\det a=\frac{1}{3!}T(a,a,a)=\frac{1}{3}(a^{\#}:a). \label{eq2/6}%
\end{equation}
If $a\in$ $H_{3}(\mathbb{O}_{\mathbb{C}})$ is written as in (\ref{eq2/1A}),
then
\begin{equation}
\det a=\alpha_{1}\alpha_{2}\alpha_{3}-\sum_{i}\alpha_{i}n(a_{i})+a_{1}%
(a_{2}a_{3})+(\widetilde{a_{3}}\widetilde{a_{2}})\widetilde{a_{1}}.
\label{eq2/7}%
\end{equation}
The following identities hold in $H_{3}(\mathbb{O}_{\mathbb{C}})$:
\begin{align}
(a^{\#})^{\#}  &  =(\det a)a,\label{eq2/8}\\
\det(a^{\#})  &  =(\det a)^{2}. \label{eq2/9}%
\end{align}

The \emph{complex conjugate }of $a=\sum_{j=1}^{3}\alpha_{j}e_{j}+\sum
_{j=1}^{3}F_{j}(a_{j})\in H_{3}(\mathbb{O}_{\mathbb{C}})$ is defined as
$\overline{a}=\sum_{j=1}^{3}\overline{\alpha_{j}}e_{j}+\sum_{j=1}^{3}%
F_{j}(\overline{a_{j}}).$ On $H_{3}(\mathbb{O}_{\mathbb{C}})$, we consider the
Hermitian scalar product defined by
\begin{equation}
\left(  a\mid b\right)  =\left(  a:\overline{b}\right)  . \label{eq2/9b}%
\end{equation}
The \emph{Jordan triple product} on $H_{3}(\mathbb{O}_{\mathbb{C}})$ is defined
by
\begin{align}
Q(x)y  &  =(x\mid y)x-x^{\#}\times\overline{y},\label{eq2/17}\\
D(x,y)z  &  =\{xyz\}=(x\mid y)z+(z\mid y)x-(x\times z)\times\overline{y}.
\label{eq2/18}%
\end{align}
With this triple product, $H_{3}(\mathbb{O}_{\mathbb{C}})$ is an Hermitian
positive Jordan triple system. The space $H_{3}(\mathbb{O}_{\mathbb{C}})$
endowed with the triple product defined by (\ref{eq2/18}) will be referred to
as the \emph{Hermitian Jordan triple system }$H_{3}(\mathbb{O}_{\mathbb{C}%
}\mathbb{)}$, or the \emph{Hermitian Jordan triple system of type VI}, or the
\emph{exceptional Hermitian triple system of dimension }$27$.

The \emph{generic minimal polynomial} of the Jordan triple system
$H_{3}(\mathbb{O}_{\mathbb{C}})$ is
\begin{equation}
m(T,x,y)=T^{3}-(x\mid y)T^{2}+(x^{\#}\mid y^{\#})T-\det x\det\overline{y};
\label{eq2/21}%
\end{equation}
the rank of $H_{3}(\mathbb{O}_{\mathbb{C}})$ is $3$. The \emph{generic norm
}is
\[
N(x,y)=1-(x\mid y)+(x^{\#}\mid y^{\#})-\det x\det\overline{y}%
\]
The \emph{set }$\mathcal{E}$\emph{ of tripotents} of $H_{3}(\mathbb{O})$ is the
disjoint
union $\mathcal{E}=\mathcal{E}_{0}\cup\mathcal{E}_{1}\cup\mathcal{E}_{2}%
\cup\mathcal{E}_{3}$, where
\begin{align*}
\mathcal{E}_{0}  &  =\{0\},\quad\mathcal{E}_{1}=\left\{  x;(x\mid
x)=1,\ x^{\#}=0\right\}  ,\\
\mathcal{E}_{2}  &  =\left\{  x;(x\mid x)=2,\ (x^{\#}\mid x^{\#})=1,\ \det
x=0\right\}  ,\\
\mathcal{E}_{3}  &  =\left\{  x;(x\mid x)=3,\ (x^{\#}\mid x^{\#}%
)=3,\ \left\vert \det x\right\vert ^{2}=1\right\}  .
\end{align*}
The elements $e_{1}$, $e_{2}$, $e_{3}$ belong to $\mathcal{E}_{1}$ and are
therefore minimal tripotents. It is easily checked that the Peirce spaces for
$e_{1}$ are
\begin{align*}
& V_{0}(e_{1})=\mathbb{C}e_{2}\oplus\mathbb{C}e_{3}\oplus\mathcal{F}_{1},\\
& V_{1}(e_{1})=\mathcal{F}_{2}\oplus\mathcal{F}_{3},\quad
V_{2}(e_{1})=\mathbb{C}e_{1}.
\end{align*}
Similar results hold for the Peirce decomposition w.r. to $e_{2}$ and $e_{3}$.
As $e_{2}$ and $e_{3}$ belong to $V_{0}(e_{1})$, they are orthogonal to
$e_{1}$; also, $e_{2}$ is orthogonal to $e_{3}$. So $(e_{1},e_{2},e_{3})$ is a
frame for the Jordan triple system $H_{3}(\mathbb{O}_{\mathbb{C}})$. It is then
easily checked that the simultaneous Peirce decomposition w.r. to the frame
$(e_{1},e_{2},e_{3})$ is
\[
H_{3}(\mathbb{O}_{\mathbb{C}})=\bigoplus_{1\leq i\leq j\leq3}V_{ij},
\]
with $V_{ii}=\mathbb{C}e_{i}$, $V_{ij}=\mathcal{F}_{k}$. As all $V_{ij}$
$(1\leq i<j\leq3)$ for this frame are non-zero, the Hermitian Jordan triple
system $H_{3}(\mathbb{O}_{\mathbb{C}})$ is \emph{simple}. Hence the numerical
invariants are
\[
a=\dim\mathcal{F}_{i}=8,\quad b=\dim V_{0i}=0,\quad r=3,\quad g=2+a(r-1)=18.
\]
In particular, $H_{3}(\mathbb{O}_{\mathbb{C}})$ is of tube type. In
$H_{3}(\mathbb{O}_{\mathbb{C}})$,
\begin{align}
\operatorname{tr}D(x,y){}  &  =18(x\mid y),\label{eq2/33}\\
\det B(x,y)){}  &  =\left(  1-(x\mid y)+(x^{\#}\mid y^{\#})-\det
x\det\overline{y}\right)  ^{18}. \label{eq2/32}%
\end{align}

The \emph{bounded circled symmetric domain} $\Omega_{VI}$ corresponding to the
Jordan triple system $H_{3}(\mathbb{O}_{\mathbb{C}})$ is defined by the set of
inequalities%
\begin{gather*}
1-(x|x)+(x^{\sharp}|x^{\sharp})-\left\vert \det x\right\vert ^{2}>0,\\
3-2(x|x)+(x^{\sharp}|x^{\sharp})>0,\\
3-(x|x)>0.
\end{gather*}

\subsection{The exceptional Jordan triple and the exceptional domain of
dimension $16$\label{Exc16}}

The $16$-dimensional vector space
\[
W=\mathcal{M}_{2,1}(\mathbb{O}_{\mathbb{C}})=\left\{  \left(  a_{2}%
,a_{3}\right)  \mid a_{2},a_{3}\in\mathbb{O}_{\mathbb{C}}\right\}
\]
is identified with the subspace of $\mathcal{H}_{3}(\mathbb{O}_{\mathbb{C}})$
consisting in matrices of the form
\[
\left(
\begin{array}
[c]{ccc}%
0 & a_{3} & \tilde{a}_{2}\\
\tilde{a}_{3} & 0 & 0\\
a_{2} & 0 & 0
\end{array}
\right)  .
\]
This is the Peirce subspace $V_{1}(e_{1})$ of $\mathcal{H}_{3}(\mathbb{O}%
_{\mathbb{C}})$ with respect to the tripotent $e_{1}$ and hence a Jordan
subsystem of the Jordan triple system $\mathcal{H}_{3}(\mathbb{O}_{\mathbb{C}%
})$. The Peirce subspace $V_{0}(e_{1})=\mathbb{C}e_{2}\oplus\mathbb{C}%
e_{3}\oplus\mathcal{F}_{1}$ will be identified with the space
\[
\mathcal{H}_{2}(\mathbb{O}_{\mathbb{C}})=\left\{  \left(
\begin{array}
[c]{cc}%
\alpha_{2} & \widetilde{a_{1}}\\
a_{1} & \alpha_{3}%
\end{array}
\right)  ;\alpha_{2},\alpha_{3}\in\mathbb{C},\ a_{1}\in\mathbb{O}_{\mathbb{C}%
}\right\}
\]
of $2\times2$ Cayley-Hermitian matrices with coefficients in $\mathbb{O}%
_{\mathbb{C}}$; the space $\mathcal{H}_{2}(\mathbb{O}_{\mathbb{C}})$ is also a
Jordan subsystem of the Jordan triple system $\mathcal{H}_{3}(\mathbb{O}%
_{\mathbb{C}})$, and it can be shown it is isomorphic to the \emph{classical}
Hermitian Jordan triple of type $IV_{10}$.

For $x=\left(  x_{2},x_{3}\right)  \in\mathcal{M}_{2,1}(\mathbb{O}%
_{\mathbb{C}})$, the adjoint $x^{\sharp}\in\mathcal{H}_{2}(\mathbb{O}%
_{\mathbb{C}})$ is then defined by
\[
x^{\sharp}=\left(
\begin{array}
[c]{cc}%
-n\left(  x_{2}\right)  & \widetilde{x_{3}}\widetilde{x_{2}}\\
x_{2}x_{3} & -n\left(  x_{3}\right)
\end{array}
\right)  .
\]
The \emph{Freudenthal product} of $x=\left(  x_{2},x_{3}\right)  ,\ y=\left(
y_{2},y_{3}\right)  \in\mathcal{M}_{2,1}(\mathbb{O}_{\mathbb{C}})$ is then
\[
x\times y=\left(  x+y\right)  ^{\sharp}-x^{\sharp}-y^{\sharp}.
\]

The spaces $\mathcal{M}_{2,1}(\mathbb{O}_{\mathbb{C}})$ and $\mathcal{H}%
_{2}(\mathbb{O}_{\mathbb{C}})$ are endowed with the following Hermitian scalar
products, inherited from the Hermitian scalar product in $\mathcal{H}%
_{3}(\mathbb{O}_{\mathbb{C}})$:

\begin{itemize}
\item if $x=\left(  x_{2},x_{3}\right)  ,\ y=\left(  y_{2},y_{3}\right)
\in\mathcal{M}_{2,1}(\mathbb{O}_{\mathbb{C}})$,
\begin{equation}
(x|y)=(x_{2}|y_{2})+(x_{3}|y_{3}); \label{EQN50}%
\end{equation}

\item if $u=\left(
\begin{array}
[c]{cc}%
\lambda_{2} & u_{1}\\
\widetilde{u_{1}} & \lambda_{3}%
\end{array}
\right)  ,\ v=\left(
\begin{array}
[c]{cc}%
\mu_{2} & v_{1}\\
\widetilde{v_{1}} & \mu_{3}%
\end{array}
\right)  \in\mathcal{H}_{2}(\mathbb{O}_{\mathbb{C}})$, ($\lambda_{2}%
,\lambda_{3},\mu_{2},\mu_{3}\in\mathbb{C}$, $u_{1},v_{1}\in\mathbb{O}%
_{\mathbb{C}}$),
\begin{equation}
(u|v)=\lambda_{2}\overline{\mu_{2}}+\lambda_{3}\overline{\mu_{3}}+(u_{1}%
|v_{1}). \label{EQN51}%
\end{equation}

\end{itemize}

The \emph{quadratic operator} of the \emph{Jordan triple system of type} $V$
is defined, for $x=\left(  x_{2},x_{3}\right)  ,\ y=\left(  y_{2}%
,y_{3}\right)  \in\mathcal{M}_{2,1}(\mathbb{O}_{\mathbb{C}})$, by
\begin{equation}
Q(x)y=\left(  x_{2}\widetilde{\overline{y_{2}}}x_{2}+(x_{2}\overline{y_{3}%
})\widetilde{x_{3}},\widetilde{x_{2}}(\overline{y_{2}}x_{3})+x_{3}%
\widetilde{\overline{y_{3}}}x_{3}\right)  . \label{EQN52}%
\end{equation}
The \emph{generic minimal polynomial} of $W$ is
\[
m(T,x,y)=T^{2}-(x\mid y)T+(x^{\#}\mid y^{\#});
\]
the \emph{generic norm }is
\[
N(x,y)=1-(x\mid y)+(x^{\#}\mid y^{\#}).
\]
The\emph{ set of tripotents} of $W$ is $\mathcal{E}^{\prime}=\mathcal{E}%
_{0}^{\prime}\cup\mathcal{E}_{1}^{\prime}\cup\mathcal{E}_{2}^{\prime}$, with
\begin{align*}
\mathcal{E}_{0}^{\prime}  &  =\{0\},\quad\mathcal{E}_{1}^{\prime}=\left\{
x\in W\mid(x\mid x)=1,\ x^{\#}=0\right\}  ,\\
\mathcal{E}_{2}^{\prime}  &  =\left\{  x\in W\mid(x\mid x)=2,\ (x^{\#}\mid
x^{\#})=1\right\}  .
\end{align*}
The triple system $W$ is simple. Its numerical invariants are%
\[
a=6,\ b=4,\ r=2,\ g=12.
\]
(The numerical invariants and the simplicity of $W$ are obtained by computing
the simultaneous Peirce decomposition w.r. to a set of two orthogonal
tripotents). In $W=\mathcal{M}_{2,1}(\mathbb{O}_{\mathbb{C}})$, we have%
\begin{align*}
\operatorname{tr}D(x,y)  &  =12(x\mid y),\\
\det B(x,y)  &  =\left(  1-(x\mid y)+(x^{\#}\mid y^{\#})\right)  ^{12}.
\end{align*}

The \emph{exceptional domain of dimension }$16$\ is
\begin{equation}
\Omega_{V}=\left\{  x\in\mathcal{M}_{2,1}(\mathbb{O}_{\mathbb{C}}%
)\mid1-(x|x)+(x^{\sharp}\mid x^{\sharp})>0,\ 2-(x|x)>0\right\}  .
\label{EQN54}%
\end{equation}

\end{document}